\documentclass[a4paper,12pt]{article}

\usepackage{amsthm,amsmath,amscd,mathtools}
\usepackage{stmaryrd}
\usepackage[vmargin=2cm,hmargin=2cm,headheight=14.5pt,top=2cm,headsep=.5cm]{geometry}
\usepackage{hyperref}
\usepackage[utf8]{inputenc}
\usepackage[autobold]{mathfixs}
\usepackage{stackrel}
\usepackage{cases}
\usepackage{circuitikz}
\usepackage{tikz}
\usepackage{tikz-cd}
\tikzcdset{every label/.append style = {font = \small}}

\usepackage[charter]{mathdesign}
\usepackage{bm}
\usepackage{caption}
\usepackage{subcaption}
\usepackage{xcolor}
\usepackage{graphicx}
\DeclareMathSizes{12}{12}{8}{6}

\usepackage{cite}
\usepackage{placeins}
\usepackage{enumitem}
\setlist[itemize,2]{label=$\centerdot$}
\setlist[itemize,3]{label=$\triangle$}

\usepackage[shortcuts]{extdash}

\newtheoremstyle{ptheorem}{1em}{0em}{\itshape}{}{\bfseries}{.}{.5em}{\thmname{#1}\thmnumber{
		#2}\thmnote{ (\hspace{-.01pt}{#3})}}

\theoremstyle{ptheorem}

\newtheorem{thm}{Theorem}[section]
\newtheorem{pro}[thm]{Proposition}
\newtheorem{lem}[thm]{Lemma}
\newtheorem{cor}[thm]{Corollary}

\newtheoremstyle{hdef}{1em}{0em}{}{}{\bfseries}{.}{.5em}{\thmname{#1}\thmnumber{
		#2}\thmnote{ (\hspace{-.01pt}{#3})}}
\theoremstyle{hdef}

\newtheorem{dfn}[thm]{Definition}
\newtheorem{rem}[thm]{Remark}

\newtheorem{exa}[thm]{Example}

\numberwithin{equation}{section}
\numberwithin{figure}{section}

\renewcommand{\phi}{\varphi}

\newcommand{\olb}[1]{%
	\vbox{\offinterlineskip\ialign{\hfil##\hfil\cr
			$\rotatebox[origin=c]{90}{$]$}$\cr\noalign{\kern-.45ex}{$#1$}\cr}}}

\newcommand{\noop}[1]{}

\usepackage{stmaryrd}

\allowdisplaybreaks

\begin{document}

\title{The Fundamental Theorem of Calculus for Lebesgue\--Stieltjes integrals involving non\--monotonic derivators}
\author{Lamiae Maia$^{1}$ and F. Adri\'an F. Tojo$^{2}$}
\date{}
\maketitle

	\begin{center} {\small $^{1}$ LAMA Laboratory, Mathematics Department, Faculty of Sciences\\
 Mohammed V University in Rabat, Rabat, Morocco. \\ e\--mail: \emph{lamiae\_maia@um5.ac.ma}}
\\
\small $^{2}$ Departamento de Estat\'{\i}stica, An\'alise Matem\'atica e Optimizaci\'on \\ Universidade de Santiago de Compostela \\ 15782, Facultade de Matem\'aticas, Campus Vida, Santiago, Spain.\\ CITMAga, Santiago de Compostela \\ e\--mail: \emph{fernandoadrian.fernandez@usc.es}
	\end{center}

	\medbreak

\begin{abstract}
In this work, we extend the concept of the Stieltjes derivative to encompass left\--continuous derivators with bounded variation, thereby relaxing the monotonicity constraint. This generalization necessitates a refined definition of the Stieltjes derivative applicable across the entire domain, accommodating derivators that may change sign. We establish a generalized Fundamental Theorem of Calculus for the Lebesgue–Stieltjes integral in this broader context, presenting both "almost\--everywhere" and "everywhere" versions. The latter requires a specific condition relating the derivator to its variation function, which we prove to be optimal through a density theorem. Our framework bridges the gap between Stieltjes differential equations and measure differential equations, offering a tool for modeling complex systems with non\--monotonic dynamics.
\end{abstract}

\medbreak


\medbreak


\section{Introduction}
In recent years, the {\em Stieltjes derivative} has attracted growing interest across the field of applied mathematics due to its versatility and ability to handle functions that are not necessarily differentiable in the classical sense. By defining the rate of change of a function~$f$ with respect to a function~$g$, the Stieltjes derivative allows for meaningful differentiation even when $g$ contains jumps or discontinuities, where standard derivatives would fail.

The idea of differentiating with respect to a function can be traced back to the early works of Daniell~\cite{D1918,D1920} and Young~\cite{Young}. This concept involves taking derivatives with respect to a function $g:\mathbb{R} \to \mathbb{R}$, known as the \emph{derivator}. Similar ideas were also explored in other works, see for instance,~\cite{AP,Cacc,Feller,Gradin,Liberman,Petrov}. Recently, in~\cite{PR}, López Pouso and Rodríguez delved into the significance of the Stieltjes derivative, coming up with a practical definition with respect to a left\--continuous and nondecreasing function. Thanks to their definition, a generalized framework unifying discrete and continuous calculus has been established, allowing to investigate solutions of ordinary differential equations, difference equations, impulsive differential equations and equations on time scales from a common standpoint~\cite{FP,ThesisFLariviere,ThesisMarquezAlbes,PR,Tojo2025connectingSDEnODE}.
 This has, in particular, led to significant applications in modeling phenomena that exhibit discontinuities and stationary periods via Stieltjes differential equations, see for instance~\cite{AFNT,FMarTo-OnFirstandSec,FernandezTojoStieltjesBochnerSpaces,FP,FrTo,PM,PM2,PM3,PMM,MEF1,MEF2,MEF3,AlbesSlavikLogisticEq}. Additionally, in~\cite{FMarTo-OnFirstandSec}, the authors refined the definition of the Stieltjes derivative with respect to a left\--continuous and nondecreasing derivator, to define it on the whole domain. Their comprehensive definition allowed considering higher\--order Stieltjes derivatives, and exploring second\--order linear Stieltjes differential with constant coefficients~\cite{FMarTo-OnFirstandSec} and non\--constant coefficients~\cite{FernándezAlbésTojo2024FirstNSec}.

In~\cite{FrTo}, the authors introduced a general setting for the Stieltjes derivative considering left\--continuous and non\--monotonic derivators, precisely those with {\em controlled variation}. Their definition takes into account the monotonic behavior of these derivators on given countable connected sets. This definition has permitted to generalize many notions from the case of monotonic derivators. However, their definition does not account for all the points of the whole domain, and the choice of derivators remains restricted to the case where the derivator exhibit monotonic behavior on those given connected components.

In the present work, we aim to explore the theory of Stieltjes derivatives beyond the monotonicity conditions typically imposed on the derivator.
Specifically, we consider left\--continuous derivators with bounded variation, which introduces additional complexities, and we provide a definition of the Stieltjes derivative across the entire domain. The motivation for choosing this level of regularity for derivators stems from the signed Lebesgue\--Stieltjes measure generated by such derivators. This would bridge the study of Stieltjes differential equations with their counterparts in measure differential equations, as in~\cite{MonteiroBianca-GenerDeriv2017}. Earlier studies of derivatives related to functions of bounded variation have been conducted from a measure\--theoretic perspective~\cite{D1918,Garg1992relativization}, as well as in connection with the Kurzweil-Stieltjes integral~\cite{MonteiroBianca-GenerDeriv2017}. The use of measures allows for a broad approach without the need to focus on the behavior at individual points. However, this same advantage can become a limitation in practical applications, where controlling the behavior at specific points is crucial via the derivative.

The Fundamental Theorem of Calculus is a basic result in analysis that relates the derivative and the integral. It is usually presented in two parts, one concerning the derivative of the integral and the other the integral of the derivative (Barrow's rule). Furthermore, the classical version of the theorem (for the Riemann integral) considers the derivative at every point of the interval of definition $[a,b]$, whereas the more modern one (the version for the Lebesgue integral) considers the derivative almost everywhere and has weaker hypotheses. In their work \cite{MonteiroBianca-GenerDeriv2017}, Monteiro and Satco have introduced the Fundamental Theorems of Calculus for Kurzweil-Stieltjes integral involving, namely, regulated functions which are $BVG^{\circ}$. This has permitted to establish an equivalence between differential problems involving distributional derivatives, Stieltjes derivatives, and those involving measure differential equations, under suitable assumptions. In our setting, we first use Lebesgue\--Stieltjes integrals and absolute continuity with respect to the integrand to establish the Fundamental Theorems of Calculus for Lebesgue\--Stieltjes integrals (the ``almost\--everywhere'' version, both for the derivative of the integral and the integral of the derivative). Second, we use our refined definition of the Stieltjes derivative to derive a new ``everywhere'' version of the Fundamental Theorem of Calculus. This version emphasizes the necessity of an assumption relating the derivator to its variation function\--a condition that is always satisfied in the case of monotonic derivators. Furthermore, by demonstrating a density theorem, we prove that this condition is optimal.

This paper is organized as follows. In the next section, we present the necessary preliminaries to define the Stieltjes derivative for non\--monotonic derivators. In Section~3, we introduce the topology and a notion of continuity related to non\--monotonic derivators. Section~4 is devoted to the concept of absolute continuity with respect to the derivator, along with related results derived via Hahn's decomposition and the Radon\--Nikod\'{y}m Theorem~\cite{BenedettoCzaja2009integration}. In Section~5, we present a generalized Fundamental Theorem of Calculus for Lebesgue\--Stieltjes integrals in both versions: the almost everywhere version and the everywhere version. The latter requires a necessary assumption, which we prove to be optimal after demonstrating a density theorem.

\section{Preliminaries}
In this preliminary section, we delve into the essential tools required to define the Stieltjes derivative with respect to $g: [a,b]\subset\mathbb{R} \to{\mathbb R}$, a left\--continuous {\em derivator} with {\em bounded variation}. In doing so, we appeal first to the Stieltjes derivative involving monotonic derivators, and some elements from measure theory related to functions with bounded variation~\cite{AL,BenedettoCzaja2009integration}.
\subsection{Monotonic derivators}
Let $[a,b]\subset{\mathbb R}$ be an interval, and $g:[a,b]\to{\mathbb R}$ a left\--continuous nondecreasing function. We will refer to such functions as \emph{derivators}. For a derivator~$g$, let $D_g=\{d_n\}_{n \in \Lambda}$, $ \Lambda\subset\mathbb{N}$, denote the set of all discontinuity points of $g$. More precisely, $D_g=\{ t \in [a,b) : \Delta^+ g(t)>0\}$ where $\Delta^+ g(t):=g(t^+)-g(t)$, $t\in\mathbb [a,b)$, and $g(t^+)$ denotes the right\--hand side limit of $g$ at $t$. Observe that $D_g$ is countable and, thus $[a,b]\backslash D_g$ is dense in $[a,b]$ (we will use this fact in the proof of Theorem~\ref{thmnecessarycondition}).

 We also define the set
\begin{equation*}
	C_g:=\{ t \in(a,b) \, : \, \mbox{$g$ is constant on $(t-\varepsilon,t+\varepsilon)$ for
		some $\varepsilon>0$} \}.
\end{equation*}
Observe that $C_g$ is open in the usual topology of $\mathbb R$, so it can be decomposed as a disjoint union of open intervals:
\begin{equation}\label{Cgdecomp}
	C_g=\bigcup_{n \in \widetilde\Lambda}(a_n,b_n),
\end{equation}
where
$\widetilde\Lambda\subset \mathbb{N}$ and $(a_k,b_k)\cap (a_j,b_j)=
\emptyset$ for $k\neq j$. With this notation, we define \[ N_g^-:=\{a_n\}_{n \in \widetilde\Lambda}\backslash (D_g\cup\{a\}),\quad N_g^+:=\{b_n\}_{n \in \widetilde\Lambda}\backslash D_g,\quad N_g:=N_g^-\cup N_g^+.\]

The derivator $g$ defines a Lebesgue\--Stieltjes measure $\mu_g$, the reader is referred to~\cite{PR,MT} for more details. We denote $\mathcal{L}^1_g([a,b),\mathbb{R})$ the space of integrable functions with respect to the measure~$\mu_g$. Moreover, we say that a property holds for {\it $g$\--almost every} $t\in I\subset [a,b)$ if it holds except on a set $N \subset I$ such that $\mu_g(N)=0$.

We now present the definition of the \emph{Stieltjes derivative}, as introduced in~\cite{FMarTo-OnFirstandSec}, which is defined on the entire domain.
 \begin{dfn}[{\cite[Definition 3.1]{FMarTo-OnFirstandSec}}]\label{dfn:Stieltjes-derivative}Let $[a,b]\subset \mathbb{R}$
	be a closed interval and $g:[a,b]\to\mathbb{R}$ a left\--continuous nondecreasing derivator such that
	$b \notin N_g^+$. We define the \emph{Stieltjes derivative}, or \emph{$g$\--derivative}, of a function $f:[a,b]\to\mathbb R$ at a point $t\in [a,b]$ as
	\begin{equation*}
		f'_g(t)=\begin{dcases}
			 \lim_{s \to t}\frac{f(s)-f(t)}{g(s)-g(t)}, &\mbox{if } t\not\in D_{g}\cup C_g,\vspace{0.1cm}\\
			\lim_{s\to t^+}\frac{f(s)-f(t)}{g(s)-g(t)}, &\mbox{if } t\in D_{g}, \vspace{0.1cm}\\
			 \lim_{s \to b_n^+}\frac{f(s)-f(b_n)}{g(s)-g(b_n)}, &\mbox{if } t \in (a_n,b_n)\subset C_g,
		\end{dcases}
	\end{equation*}
with $a_n,b_n$ as in~\eqref{Cgdecomp}; provided the corresponding limits exist. In that case, we say that $f$ is \emph{$g$\--differentiable at~$t$}.
Furthermore, the $g$\--derivative at a point $t\in N_g\cup\{a,b\}$ must be understood as
\begin{equation*}
		f'_g(t)=\begin{dcases}
			\lim_{\substack{s \to t^+\\ g(s)\ne g(t)}}\frac{f(s)-f(t)}{g(s)-g(t)}, &\mbox{if } t\in N_g^+\cup\{a\},\\
			\lim_{s\to t^-}\frac{f(s)-f(t)}{g(s)-g(t)}, &\mbox{if } t\in N_g^-\cup\{b\}.
		\end{dcases}
\end{equation*}
\end{dfn}
Observe that, in the first case, we are including the case $a\notin D_g$ but $(a,b_n)\subset C_g$. In the sequel, for all $t\in[a,b]$, we consider the following notation introduced in~\cite[Proposition~3.9]{FMarTo-OnFirstandSec}:
\begin{equation}\label{eq:notation of t^*}
 t^*=\begin{dcases}
 b_n, & \mbox{if $t\in(a_n,b_n)\subset C_g$, or ($t=a\notin D_g$ and $(a,b_n) \subset C_g$)},\\
 t, & \mbox{otherwise.}
 \end{dcases}
\end{equation}
Observe that $t^*\in [a,b]\setminus C_g$ for all $t\in[a,b]$. Using this notation, we can express the definition of the Stieltjes derivative for all $t\in(a,b)$ as
	\begin{equation*}
		f'_g(t)=\begin{dcases}
			 \lim_{s \to t^*}\frac{f(s)-f(t^*)}{g(s)-g(t^*)}, &\mbox{if } t^*\not\in (D_g\cup N_g),\\
			\lim_{s\to t^{*+}}\frac{f(s)-f(t^*)}{g(s)-g(t^*)}, &\mbox{if } t^*\in (D_g\cup N_g^+),\\
			\lim_{s\to t^{*-}}\frac{f(s)-f(t^*)}{g(s)-g(t^*)}, &\mbox{if } t^*\in N_g^-.
		\end{dcases}
	\end{equation*}
Now we characterize Definition~\ref{dfn:Stieltjes-derivative} with the following proposition.
\begin{pro}\label{prop:Stieltjes-differentiability-monotonic-derivator}
 Let $[a,b]\subset \mathbb{R}$ be a closed interval and $g:[a,b]\to\mathbb{R}$ a left\--continuous nondecreasing derivator such that $b\notin N_g^+$ and let $t\in [a,b]$. Then the following statements are equivalent
\begin{enumerate}
 \item $f:[a,b]\to{\mathbb R}$ is $g$\--differentiable at $t$.
 \item There exist $d\in{\mathbb R}$ and a function $h:[a,b]\to {\mathbb R}$ satisfying
\begin{enumerate}
 \item $h(t^*)=0$;
 \item $h$ is continuous at $t^*$ if $t^*\in [a,b]\setminus(D_g\cup N_g)$;
 \item $h$ is right\--continuous at $t^*$ if $t^*\in(D_g\cup N_g^+)$;
 \item $h$ is left\--continuous at $t^*$ if $t^*\in N_g^-$;
\end{enumerate}
\end{enumerate}
 such that, $f(s)=f(t^*)+[d+h(s)][g(s)-g(t^*)]$ for $s\in[a,b]$, $g(s)\ne g(t^*)$,		with $t^*$ as in~\eqref{eq:notation of t^*}. When these properties hold, $d=f'_g(t)$.
\end{pro}
\begin{proof} $1\Rightarrow 2$. Let $t\in [a,b]$, and assume that $f$ is $g$\--differentiable at $t$. Let us consider the function $h:[a,b]\to{\mathbb R}$ defined by
\[
h(s)=\begin{dcases}
 \frac{f(s)-f(t^*)}{g(s)-g(t^*)}-f_g'(t), & \mbox{if } g(s)\neq g(t^*), \\
 0, & \mbox{otherwise}.
 \end{dcases}
\]
$h$ is well\--defined.

 $\diamond$ If $t^*\in [a,b]\setminus(D_g\cup N_g)$, then $t^*=t$ and
\[
f_g'(t)=\lim_{s\to t^*}\frac{f(s)-f(t^*)}{g(s)-g(t^*)} \in{\mathbb R}.
\]
Thus, for every $\varepsilon>0$, there exists $\delta>0$ such that
\[
\left|\frac{f(s)-f(t^*)}{g(s)-g(t^*)}-f_g'(t)\right|<\varepsilon,
\]
for $0<|s-t^*|<\delta$. Therefore, for every $\varepsilon>0$, there exists $\delta>0$ such that
\[
|h(s)-h(t^*)|<\varepsilon,
\]
for $0<|s-t^*|<\delta$. Hence, $h$ is continuous at $t^*$.

 $\diamond$ If $t^*\in D_g\cup N_g^+$, then
\[
f_g'(t)=\lim_{s\to t^{*+}}\frac{f(s)-f(t^*)}{g(s)-g(t^*)}\in{\mathbb R}.
\]
Thus, for every $\varepsilon>0$, there exists $\delta>0$ such that
\[
\left|\frac{f(s)-f(t^*)}{g(s)-g(t^*)}-f_g'(t)\right|<\varepsilon,
\]
for $0<s-t^*<\delta$. Therefore, for every $\varepsilon>0$, there exists $\delta>0$ such that
\[
|h(s)-h(t^*)|<\varepsilon,
\]
for $0<s-t^*<\delta$. Hence, $h$ is right\--continuous at $t^*$.

 $\diamond$ If $t^*\in N_g^-$, then arguing as in the second point, we deduce that $h$ is left\--continuous at $t^*$.

$2\Rightarrow 1$.
Let us assume that there exist $d\in{\mathbb R}$, and $h:[a,b]\to{\mathbb R}$ satisfying Conditions~(a)\--(d) such that, $f(s)=f(t^*)+[d+h(s)][g(s)-g(t^*)]$ for $s\in[a,b]$, $g(s)\ne g(t^*)$, that is,
\[
h(s)= \frac{f(s)-f(t^*)}{g(s)-g(t^*)}-d, \mbox{if } g(s)\neq g(t^*).
\]
Let $t\in[a,b]$. If $t^*\notin D_g\cup N_g$ then $t^*=t$. Since $h$ is continuous at $t^*$ and $h(t^*)=0$, then
\[
\lim_{s \to t^*}\frac{f(s)-f(t^*)}{g(s)-g(t^*)}=\lim_{s \to t^*} h(s)+d=d.
\]
Now, for $t^*\in D_g\cup N_g^+$ (resp. $t^*\in N_g^-$), then by the right\--continuity (resp. left\--continuity) of~$h$ at~$t^*$, we obtain the analogous result using the right\--hand limit (resp. the left\--hand limit). It follows from each case that there exists $d=f_g'(t) \in{\mathbb R}$. Hence $f$ is $g$\--differentiable at $t$.
\end{proof}

Later on, we will use the following versions of the Fundamental Theorem of Calculus, the first concerning the Stieltjes derivative of the integral and the second the integral of the Stieltjes derivative.
\begin{thm}[{\cite[Theorem 2.4]{PR}}]
	\label{derintstieljtes}
Let $[a,b]\subset \mathbb{R}$ be a closed interval, $g:[a,b]\to\mathbb{R}$ a left\--continuous nondecreasing derivator, and $f\in\mathcal{L}^1_g([a,b),\mathbb{R})$. Consider the function $F:[a,b]\to \mathbb{R}$ given by
	\[  F: x\in[a,b]\to F(x)=\int_{[a,x)}f\operatorname{d}\mu_g.\]
	Then, there exists $N\subset[a,b)$ such that $\mu_g(N)=0$ and $F'_g(x)=f(x)$ for all $x\in[a,b)\setminus N$.
\end{thm}
\begin{thm}[{Fundamental Theorem of Calculus for the Lebesgue\--Stieltjes Integral \cite[Theorem~5.4]{PR}}]\label{ftc}
Let $[a,b]\subset \mathbb{R}$ be a closed interval, $g:[a,b]\to\mathbb{R}$ a left\--continuous nondecreasing derivator, and let $F:[a,b] \to \mathbb{R}$. The following statements are equivalent:
\begin{enumerate}
 \item The function $F$ is \emph{$g$\--absolutely continuous} on $[a, b]$, i.e. for each $\varepsilon > 0$, there exists $\delta > 0$ such that, for any family $\{(a_i , b_i)\}_{i=1}^{i=m}$ of pairwise disjoint open subintervals of $[a,b]$,
\[
\sum_{i=1}^{m} g(b_i)-g(a_i) <\delta \Rightarrow \sum_{i=1}^{m} |F(b_i)-F(a_i)| < \varepsilon.
\]
 \item The following three statements are fulfilled:
	\begin{enumerate}\item There exists $F_g'(t)$ for $g$\--almost every $t \in[a, b)$;
		\item $F_g' \in \mathcal{L}_g^1([a, b),\mathbb{R})$;
		\item For each $t \in[a, b]$, we have
		\[
		F(t)=F(a)+\int_{[a, t)} F_g' \operatorname{d} \mu_g .
		\]
	\end{enumerate}
\end{enumerate}
\end{thm}
\subsection{General non\--monotonic derivators of bounded variation}
In this subsection, we consider non\--monotonic derivators of bounded variation, and explore implications of the Hahn's and Jordan's decomposition\---key tools for our study of the Generalized Fundamental Theorem of Calculus in Section~5. Let $[a,b]$ be a fixed interval of ${\mathbb R}$.
\begin{dfn}
Let $g:[a,b]\to {\mathbb R}$ be a function. Consider ${\mathcal P}([a,b])$ the set of the partitions of the interval $[a,b]$, i.e.
\[
{\mathcal P}([a,b]):=\{ P=(t_1,\dots,t_{n_P}):\, n_P\geqslant 2;\, t_i\in [a,b],\, i=1,\dots, n_P\, ; a=t_1\leqslant\dots\leqslant t_{n_P}=b\}.
\]
We define the {\em total variation} of $g$ in $[a,b]$ by
\[
\operatorname{var}_g [a,b]=\sup_{P\in {\mathcal P}([a,b])}\sum_{i=1}^{n_P -1}|g(t_i)-g(t_{i+1})|.
\]
If $\operatorname{var}_g [a,b] < \infty$, then $g$ is said to be a function of {\em bounded variation}, and we denote by $\operatorname{BV}([a,b],{\mathbb R})$ the set of functions of bounded variation and by $\operatorname{\operatorname{BV}^{lc}}([a,b],{\mathbb R})$ the set of those in $\operatorname{BV}([a,b],{\mathbb R})$ that are left\--continuous.
\end{dfn}
\begin{rem}
 Clearly, if $g\in \operatorname{BV}([a,b],{\mathbb R})$ is a nondecreasing function, then
 \[
 \operatorname{var}_g [a,b]=g(b)-g(a).
 \]
\end{rem}

\begin{dfn}\label{dfn:variation function}
 Let $g\in\operatorname{BV}([a,b],{\mathbb R})$. We define the {\em variation function} of $g$ by the function $\widetilde{g}:[a,b]\to{\mathbb R}$ given as
 \[
 \widetilde{g}(t):=\operatorname{var}_g [a,t].
 \]
\end{dfn}
\begin{rem}\ \label{rem:Dg+(t)=Dtil(g)+(t) + var[a,b]=til(g)(b)-tilde(a)}
\begin{itemize}
 \item The function $\widetilde{g}$ is nondecreasing. Furthermore, $\widetilde{g}$~is constant on the intervals where~$g$ is constant and, at any discontinuity point~$t$ of~$g$, the variation function~$\widetilde{g}$ will experience a discrete jump. In addition, for $\alpha,\beta\in[a,b]$, with $\alpha<\beta$, it is easy to verify that
\begin{equation}\label{eq:var[a,b]=tilde(g)(b)-tilde(g)(a)}
 \operatorname{var}_g [\alpha,\beta]=\widetilde{g}(\beta)-\widetilde{g}(\alpha).
\end{equation}
 \item If $g$ is left\--continuous, then so is~$\widetilde{g}$. In this case, at any point $t\in D_g$, $|\Delta^+g(t)|=\Delta^+\widetilde{g}(t)$. To not complicate the notation, since $\widetilde{g}$~shares the same discontinuity points of~$g$ and its constancy intervals, we set
\begin{equation*}
 D_g:=D_{\widetilde{g}},\quad C_g:=C_{\widetilde{g}},\quad N_g^{\pm}:=N_{\widetilde{g}}^{\pm},\quad\text{and}\quad N_g:=N_{g}^+\cup N_g^-.
\end{equation*}
\end{itemize}
\end{rem}

In the next theorem, we recall the Jordan decomposition theorem, the reader is referred to\cite[Theorem~4.1.2]{BenedettoCzaja2009integration} for further details.
\begin{thm}[Jordan decomposition theorem]\label{jdt}
If $F\in\operatorname{BV}^{\operatorname{lc}}([a,b],{\mathbb R})$, then there exist left\--continuous nondecreasing functions $F_1,F_2: [a,b] \to{\mathbb R}$ such that $F=F_1-F_2$.
\end{thm}

In the sequel, we consider a derivator $g\in\operatorname{BV}^{\operatorname{lc}}([a,b],{\mathbb R})$. The derivator $g$, being a function of bounded variation, generates a signed measure~$\mu_g$ on the measurable space $([a,b],\mathcal{M}_g)$ where $\mathcal{M}_g$ is the Borel $\sigma$\--algebra induced by the usual topology on $[a,b]\subset\mathbb{R}$. The definition of $\mu_g$ starts with the fundamental formula
\begin{equation}\label{eq:mu_g([a,t))=g(t)}
 \mu_g([a,t))=g(t)-g(a),\quad\text{for all }t\in [a,b].
\end{equation}
We refer to the measurability with respect to the measure $\mu_g$ by $g$\--measurability.
The Jordan decomposition will play a key role in decomposing the derivator $g\in \operatorname{BV}^{\operatorname{lc}}([a,b],{\mathbb R})$ into a difference of two monotone nondecreasing derivators $g_1,g_2:[a,b]\to{\mathbb R}$ and relate the Lebesgue\--Stieltjes measure associated to each derivator with Hahn's decomposition~\cite[Theorem~5.1.6]{BenedettoCzaja2009integration}. We summarize this result in the following theorem; the reader is referred to~\cite[Theorems~5.1.6,~5.1.8, and~5.1.9]{BenedettoCzaja2009integration} for further details.

\begin{thm}[Consequence of the Hahn's and Jordan's decomposition]\label{thm:consequence:Hahn+Jordan Decomp-sets A+g & A-g}
Consider the measure space $([a,b],{\mathcal M}_g,\mu_g)$. Then, there exist $g$\--measurable sets $A_g^+$, $A_g^-$ such that:
\begin{enumerate}
 \item $A_g^+ \cap A_g^- = \emptyset$, and $A_g^+\cup A_g^-=[a,b]$;
 \item $\mu_g(E) \geqslant 0$ if $E\subset A_g^+$ is $g$\--measurable;
 \item $\mu_g(E) \leqslant 0$ if $E\subset A_g^-$ is $g$\--measurable.
\end{enumerate}
 If we set
 \[  \mu_g^+(E):=\mu_g(E\cap A_g^+),\text{ and } \mu_g^-(E):=-\mu_g(E\cap A_g^-)\quad \text{for all }E\in {\mathcal M}_g,
 \]
 then $\mu_g^+$ and $\mu_g^-$ are positive measures on ${\mathcal M}_g$ and $ \mu_g=\mu_g^+ - \mu_g^-$.
 The measure $\mu_g^+$ (resp. $\mu_g^-$) is called the \emph{positive} (resp. {\em negative}) variation of the measure $\mu_g$. The measure $|\mu_g|$ defined by $
|\mu_g|:=\mu_g^+ + \mu_g^-
$,
is called the {\em total variation} of $\mu_g$.
\end{thm}

Now, given Remark~\ref{rem:Dg+(t)=Dtil(g)+(t) + var[a,b]=til(g)(b)-tilde(a)}, formula~\eqref{eq:mu_g([a,t))=g(t)}, and Theorem~\ref{thm:consequence:Hahn+Jordan Decomp-sets A+g & A-g}, it is immediate to establish the following corollary.
\begin{cor}\label{cor:relating mu_g+ and mu_g1, mu_g- and mu_g2}
Let $g\in \operatorname{BV}^{\operatorname{lc}}([a,b],{\mathbb R})$. Let us define the functions $g_1, g_2:[a,b]\to \mathbb{R}$ by
\[
g_1(t):=\mu_g^+([a,t))\quad\text{and}\quad g_2(t):=\mu_g^-([a,t)),\quad \text{for all }t\in [a,b].
\]
The functions $g_1$ and $g_2$ are left\--continuous and nondecreasing derivators. Moreover, $g=g_1 -g_2 +g(a)$ and the following statements hold:
\begin{enumerate}
 \item $\mu_{g_1}=\mu_g^+$, and $\mu_{g_2}=\mu_g^-$, where $\mu_{g_1}$ (resp. $\mu_{g_2}$) is the Lebesgue\--Stieltjes measure generated by the derivator $g_1$ (resp. $g_2$).
 \item $\mu_g^\pm\left( A_g^\mp \right) =0$.
 \item The measure $|\mu_g|$ coincides with the Lebesgue\--Stieltjes measure $\mu_{\widetilde{g}}$ generated by the derivator~$\widetilde{g}$, and for all $x,y\in[a,b]$, with $x<y$, $|\mu_g|$ satisfies
 \[
 |\mu_g|([x,y))=\widetilde{g}(y)-\widetilde{g}(x)=\operatorname{var}_g[x,y].
 \]
\end{enumerate}
\end{cor}
We refer to the \emph{measurability} (resp. \emph{integrability}) with respect to the measure $|\mu_g|$ by $\widetilde{g}$\--measurability (resp. $g$\--integrability), and we denote \ $L_g^1(I,\mathbb{R})$ the space of $g$\--integrable functions on the interval $I=[a,b)\subset \mathbb{R}$ endowed with the norm
\[
\|f\|_{L_{g}^{1}(I)}:=\int_{I} |f| \operatorname{d} |\mu_g|, \quad \text{ for every } f\in L_{g}^{1}(I,\mathbb{R}).
\]
We say that a property holds for \emph{$|\mu_g|$\--almost every $t\in I$} or \emph{$|\mu_g|$\--almost everywhere} (shortly, $|\mu_g|$\--a.e.) if it holds except on a set $N\subset I$ such that $|\mu_g|(N)=0$.
\subsection{The Stieltjes derivative for non\--monotonic derivators}
In this subsection, we generalize the definition of the Stieltjes derivative to non\--monotonic derivators of bounded variation.
Throughout this part, let $[a,b]\subset \mathbb{R}$ be a closed interval and $g:[a,b] \to \mathbb{R}$ a derivator of $\operatorname{BV}^{\operatorname{lc}}([a,b],{\mathbb R})$ such that $b\notin N_g^+$.
\begin{dfn}\label{dfn:g-derivative(non monotonic case)}
Let $f:[a,b]\to\mathbb R$, $t\in[a,b]$ and assume there exist $d\in{\mathbb R}$, and a function $h:[a,b]\to \mathbb{R}$ satisfying
\begin{enumerate}
 \item $h(t^*)=0$;
 \item $h$ is continuous at $t^*$ if $t^*\in [a,b]\setminus(D_g\cup N_g)$;
 \item $h$ is right\--continuous at $t^*$ if $t^*\in D_g\cup N_g^+$;
 \item $h$ is left\--continuous at $t^*$ if $t^*\in N_g^-$;
\end{enumerate}
 such that, $f(s)=f(t^*)+[d+h(s)][g(s)-g(t^*)]$ for $s\in[a,b]$, $g(s)\ne g(t^*)$
 with $t^*$ as in~\eqref{eq:notation of t^*}. In that case we say that $f'_g(t)\equiv d$ is the {\em Stieltjes derivative} or {\em$g$\--derivative} of $f$ at $t$ and that $f$ is \emph{$g$\--differentiable at $t$}.
\end{dfn}
\begin{rem}\label{rem:ld}Observe that Definition~\ref{dfn:g-derivative(non monotonic case)} is equivalent to defining
	\[  f'_g(t):=\begin{dcases}
		\lim_{\substack{s \to t^*\\ g(s)\ne g(t^*)}}\frac{f(s)-f({t^*})}{g(s)-g({t^*})}, &\mbox{if } {t^*}\not\in D_{g},\\
		\lim_{s\to {t^{*+}}}\frac{f(s)-f({t^*})}{g(s)-g({t^*})}, &\mbox{if } {t^*}\in D_{g},\end{dcases}\]
cf. \cite[Definition 2.26]{MonteiroBianca-GenerDeriv2017}. In particular, in the case where $t^*\in D_g$ in Definition~\ref{dfn:g-derivative(non monotonic case)}, since $g$ is regulated then $f'_g(t)$ exists if and only if $f(t^{*+})$ exists, and we have
\[
f'_g(t)=\frac{f(t^{*+})-f(t^*)}{\Delta^+ g(t^*)}.
\]
\end{rem}

\begin{rem}
 The derivator $g$ is not necessarily $\widetilde{g}$\--differentiable, and conversely, $\widetilde{g}$ is not necessarily differentiable with respect to~$g$, indicating that the Stieltjes differentiability with respect to one does not imply the Stieltjes differentiability with respect to the other. For instance for $g:{[0,2]}\to\mathbb{R}$ defined by
\begin{equation}\label{eq:rem:topologies comparison: g}
g(t)=\begin{dcases}
 t & \mbox{if } t\leqslant 1, \\
 2-t & \mbox{otherwise}.
 \end{dcases}
\end{equation}
Notice that $\widetilde{g}={\operatorname{id}_{[0,2]}}$ and $g$ is not $\widetilde{g}$\--differentiable at~$1$. Moreover,
\[
\lim\limits_{t\to 1^-} \frac{\widetilde{g}(t)-\widetilde{g}(1)}{g(t)-g(1)}=1\neq-1=\lim\limits_{t\to 1^+} \frac{\widetilde{g}(t)-\widetilde{g}(1)}{g(t)-g(1)}.
\]
Thus, $\widetilde{g}$ is not $g$\--differentiable at~$1$.
\end{rem}

\section{$g$\--topology and $g$\--continuity}
This subsection is devoted to introduce the topology generated by a derivator $g\in \operatorname{BV}^{\operatorname{lc}}([a,b],{\mathbb R})$, as well as an interesting continuity notion that can be leveraged in deriving the generalized fundamental theorem of calculus, particularly in relation to its everywhere version.

Given a derivator $g\in{\operatorname{BV}^{\operatorname{lc}}([a,b],{\mathbb R})}$, $g$ defines a pseudometric $\rho_g:{[a,b]\times[a,b]} \to \mathbb{R}^+$ given, for $s,t \in{[a,b]}$, by
\[
\rho_g(s,t)=|\widetilde{g}(s)-\widetilde{g}(t)|=\operatorname{var}_g [\min\{t,s\},\max\{t,s\}].
\]
Observe that $\rho_g(s,t)=|\Delta_g(s,t)|$, where $\Delta_g(s,t)=\widetilde{g}(s)-\widetilde{g}(t)$ is a \emph{displacement}~\cite[Definition~2.12]{ThesisMarquezAlbes}. Thus, the pseudometric $\rho_g$ generates a topology, which we denote $\tau_g$, over {$[a,b]$} given by its basic neighborhoods at each point $t\in{[a,b]}$ by the $g$\--open balls
\[
B_g(t,r)=\{s\in[a,b]: \rho_g (s,t)<r\}.
\]

\begin{rem}
 The topology $\tau_g$ is Hausdorff if and only if $\widetilde{g}$ is injective. In this case, $C_g=\emptyset$.
\end{rem}

In the following definition, we define the notion of $g$\--continuity with respect to $g$.
\begin{dfn}
Let $I\subset{[a,b]}$. A function $f:I \to {\mathbb{{\mathbb R}}}$ is said to be {\em $g$\--continuous} at $t \in I$, if, for every $\varepsilon>0$, there exists $\delta>0$ such that for all $s\in I$,
 \[
 s\in B_g(t,\delta) \Rightarrow |f(s)-f(t)|<\varepsilon.
 \]
We denote by $\mathcal{C}_{g}([a,b],{\mathbb R})$ the space of $g$\--continuous functions, and by $\mathcal{BC}_{g}([a,b],{\mathbb R})$ the space of $g$\--continuous functions which are bounded on the interval $[a,b]$. Analogously to~\cite[Theorem~3.4]{FP}, the space $\mathcal{BC}_{g}([a,b],{\mathbb R})$ equipped with supremum norm
\[
\|f\|_0=\sup_{t\in[a,b]}|f(t)|,\quad\text{for all } f\in\mathcal{BC}_{g}([a,b],{\mathbb R}),
\]
is a Banach space.
\end{dfn}

We recall the following proposition, which, although originally stated for the case of a nondecreasing $g$, is valid in the general case.
\begin{pro}[{\cite[Proposition~3.2]{FP}}]\label{prop:g-cont fonc properties}
{Let} $f:[a,b]\to \mathbb{R}$ be a $g$\--continuous function on $[a,b]$. Then the following statements hold:
\begin{enumerate}
 \item $f$ is continuous from the left at each $t\in (a,b]$.
 \item If $g$ is continuous at $t\in [a,b]$, then so is $f$.
 \item If $g$ is constant on some interval $[u,v]\subset [a,b]$, then so is $f$.
\end{enumerate}
\end{pro}

\begin{rem}\label{rem:measurability of g-cont functs}
It is worth noting that a $g$\--continuous function $f:I\to \mathbb{R}$ defined on a Borel set $I\subset {[a,b]}$ is $g$\--measurable since, in this case, $f:(I,\tau_g) \to(\mathbb{R},\tau_u)$ is continuous, so it follows that $f$ is Borel measurable. Using the same argument as in~\cite[Corollary~3.5]{FP}, we conclude that $f$ is Lebesgue\--Stieltjes measurable.
\end{rem}

While the primary focus of this section is on $g$\--continuity defined via the variation function $\widetilde{g}:{[a,b]}\to\mathbb{R}$, which involves the $g$\--topology derived from the pseudometric~$\rho_g$, it is worth noting that one could consider employing an alternative topology generated by~$g$ to define $g$\--continuity: the topology induced by pseudometric~$\breve{\rho}_g:[a,b]\times[a,b]\to \mathbb{R}^+$ given, for $s,t \in{[a,b]}$, by
\[
\breve{\rho}_g(s,t)=|g(s)-g(t)|,
\]
and related to the displacement $\Delta(t,s):=|g(s)-g(t)|$. Thus, $\breve{\rho}_g$ defines a topology $\breve{\tau}_g$ over $[a,b]$ with a local open neighborhood basis at each point $t\in {[a,b]}$ given by the $\breve{\tau}_g$\--open balls
\[
\breve{B}(t,r)=\{ {s\in[a,b]}: \breve{\rho}_g (s,t)<r\}.
\]
However, this alternative definition raises concerns regarding its applicability. Specifically, if we define a function $f:I\subset {[a,b]} \to \mathbb{R}$ at $t\in I$ to be $\breve{g}$\--continuous if, for every $\varepsilon>0$, there exists $\delta>0$ such that for all $s\in I$,
 \[
 s\in \breve{B}(t,\delta) \Rightarrow |f(s)-f(t)|<\varepsilon,
 \]
we are restricting set of $g$\--continuous functions to a smaller subset. The implications of this restriction are considered in the next remark.

\begin{rem}\label{rem:g-cont implies tilde(g)-cont}\
\begin{itemize}
\item The $\breve{g}$\--continuity defined above implies $g$\--continuity. Indeed, assume that $f:[a,b]\to\mathbb{R}$ is $\breve{g}$\--continuous at $t\in[a,b]$. Then for every $\varepsilon>0$ there exists $\eta>0$ such that
\[
|g(s)-g(t)|<\eta \Rightarrow |f(s)-f(t)|<\varepsilon.
\]
Now, using the inequality
\[
|g(s)-g(t)| \leqslant |\widetilde g(s)-\widetilde g(t)| = \rho_g(s,t), \quad \text{for every $s,t\in[a,b]$,}
\]
we obtain that for every $s,t\in[a,b]$,
\[
\rho_g(s,t)<\eta \Rightarrow |g(s)-g(t)|<\eta.
\]
Combining the two implications yields
\[
\rho_g(s,t)<\eta \Rightarrow |f(s)-f(t)|<\varepsilon,
\]
which proves that $f$ is $g$\--continuous at $t$.
 \item The converse does not necessarily hold, that is, $g$\--continuous functions may not be $\breve{g}$\--continuous. Indeed, reconsider the derivator~$g$ defined in~\eqref{eq:rem:topologies comparison: g} where $\widetilde{g}=\operatorname{id}_{{[0,2]}}$. $\widetilde{g}$ is clearly $g$\--continuous, however it is not $\breve{g}$\--continuous. To see this, observe that $\breve{g}$\--continuity can be regarded as the continuity between the topological spaces $([0,2],\breve{\tau}_g)$ and $(\mathbb{R},\tau_u)$ where~$\tau_u$ is the usual topology in $\mathbb{R}$. Thus, if $\widetilde{g}$ were indeed $\breve{g}$\--continuous, it would imply that $\widetilde{g}:([0,2],\breve{\tau}_g)\to(\mathbb{R},\tau_u)$ is continuous. However, observe that $\breve{\rho}_g(0,2)=|g(0)-g(2)|=0$, so~$0$ and~$2$ belong to the same $\breve{\rho}_g$\--equivalence class induced by the pseudometric, i.e. $0\sim 2$. Since $(\mathbb{R},\tau_u)$ is Hausdorff, continuity of~$\widetilde{g}$ would imply that $\widetilde{g}(0)=\widetilde{g}(2)$. However, $\widetilde{g}(0)=0$ and $\widetilde{g}(2)=2$, which is a contradiction. Hence, $\widetilde{g}$ is not $\breve{g}$\--continuous.
\end{itemize}
\end{rem}

\section{$g$\--absolute continuity}
In the classical context of differentiation, some functions can be reconstructed by integrating their derivatives, a process that is feasible only if the function is absolutely continuous. In the framework of the Stieltjes derivative, a similar concept was established via the idea of absolute continuity with respect to a left\--continuous and nondecreasing function $g : \mathbb{R} \to \mathbb{R}$, as defined in~\cite[Definition~5.1]{PR}. In the context of non\--monotonic derivators of {bounded} variation, a similar notion was introduced in~\cite[Definition~6.1]{FrTo} considering particular derivators $g : I \to \mathbb{R}$ with {\em controlled variation}~\cite[Definition~3.1]{FrTo} on $I\subset\mathbb{R}$. This definition takes into account the monotonic behavior of the derivator on given countable connected sets of~$I$. In the following definition, we { fix $[a,b]\subset\mathbb{R}$, and we continue to consider a derivator $g\in \operatorname{BV}^{\operatorname{lc}}([a,b],{\mathbb R})$} which is not necessarily of {\em controlled variation}.
\begin{dfn}
A map $F:[a,b]\to \mathbb{R}$ is {\em $g$\--absolutely continuous} if, for every $\varepsilon > 0$, there exists $\delta > 0$ such that, for any family
$\{(a_i , b_i)\}_{i=1}^{i=m}$ of pairwise disjoint open subintervals of $[a,b]$,
\[
\sum_{i=1}^{m} \operatorname{var}_g[a_i,b_i] <\delta \Rightarrow \sum_{i=1}^{m} |F(b_i)-F(a_i)| < \varepsilon.
\]
\end{dfn}
We denote by $\mathcal{AC}_{g}([a,b],{\mathbb R})$ the set of $g$\--absolutely continuous on the interval $[a,b]$.

\begin{rem}\label{rem:g-abs-cont is tilde(g)-abs-cont}
 In light of equation~\eqref{eq:var[a,b]=tilde(g)(b)-tilde(g)(a)}, it is worthwhile to mention that our definition of $g$\--absolute continuity coincides with $\widetilde{g}$\--absolute continuity in the sense of~\cite[Definition~5.1]{PR}. This fact, combined with~\cite[Proposition~5.5]{FP}, implies that $\mathcal{AC}_{g}([a,b],{\mathbb R})\subset\mathcal{BC}_{g}([a,b],{\mathbb R})$.
\end{rem}

In particular,~\cite[Proposition~5.3]{PR} yields the following lemma.
\begin{lem}\label{lem:g-abs-cont func has bounded var}
 $\mathcal{AC}_{g}([a,b],{\mathbb R}) \subset\operatorname{\operatorname{BV}^{lc}}([a,b],{\mathbb R})$.
\end{lem}

\begin{thm}\label{thm: F is ACg is equiv to mu_F satisfying 2 properties} Let $F\in \operatorname{\operatorname{BV}^{lc}}([a,b],{\mathbb R})$, and $\mu_F$ be the Borel measure generated by the function~$F$ as in~\eqref{eq:mu_g([a,t))=g(t)}.
 The following statements are equivalent:
\begin{enumerate}
\item $F$ is $g$\--absolutely continuous on $[a,b]$;
\item The measure $\mu_F$ satisfies {for every $\widetilde{g}$\--measurable set $E$:}
 \begin{enumerate}
 \item If $\mu_g^+(E)=0$, then $\mu_F(E\cap A_g^+)=0$;
 \item If $\mu_g^-(E)=0$, then $\mu_F(E\cap A_g^-)=0$.
 \end{enumerate}
\end{enumerate}
\end{thm}

\begin{proof}
 $(1)\Rightarrow(2)$. If $F$ is $g$\--absolutely continuous on $[a,b]$, then, by Remark~\ref{rem:g-abs-cont is tilde(g)-abs-cont}, $F$ is $\widetilde{g}$\--absolutely continuous on $[a,b]$ and $\mu_F \ll \mu_{\widetilde{g}}$. Thus, by means of the Radon\--Nikod\'{y}m~\cite[Theorem~5.3.2]{BenedettoCzaja2009integration}, there exists $h\in L^{1}_{\widetilde{g}}([a,b),\mathbb{R})$ such that $\mu_ F(E)=\int_E h\,\operatorname{d} \mu_{\widetilde{g}}=\int_E h\,\operatorname{d} |\mu_g|$ for all $\widetilde{g}$\--measurable set $E$. Thus, for {every} $\widetilde{g}$\--measurable set $E$, we have
\[
\mu_F(E)= \int_E h\, \operatorname{d} \mu_g^+ + \int_E h\, \operatorname{d} \mu_g^-.
\]
By Corollary~\ref{cor:relating mu_g+ and mu_g1, mu_g- and mu_g2}(2), $\mu_g^-(A_g^+)=0$ and $\mu_g^+(A_g^-)=0$. Therefore,
\[
\mu_F(E\cap A_g^+)= {\int_{E\cap A_g^+} h\, \operatorname{d} \mu_g^+=}\int_E h\, \operatorname{d} \mu_g^+
\text{ and }
\mu_F(E\cap A_g^-)= {\int_{E\cap A_g^-} h\, \operatorname{d} \mu_g^-=}\int_E h\, \operatorname{d} \mu_g^-.
\]
Hence,
\[
\mu_g^+(E)=0 \Rightarrow\mu_F(E\cap A_g^+)=0,
\]
and
\[
\mu_g^-(E)=0 \Rightarrow\mu_F(E\cap A_g^-)=0.
\]
$(2)\Rightarrow(1)$. Assume that $(2)$ holds. Let $E$ be a $\widetilde{g}$\--measurable set. If $|\mu_g|(E)=0$, then $\mu_g^+(E)=0$ and $\mu_g^-(E)=0$. Thus, $\mu_F(E\cap A_g^+)=0$ and $\mu_F(E\cap A_g^-)=0$. Therefore, $\mu_F(E)=0$. Hence, $\mu_F\ll \mu_{\widetilde{g}}$ and $F$ is $\widetilde{g}$\--absolutely continuous on $[a,b]$, which is $g$\--absolutely continuous on $[a,b]$ by Remark~\ref{rem:g-abs-cont is tilde(g)-abs-cont}.
\end{proof}
Consequently, we obtain the following corollary.
\begin{cor}\label{cor:f g-abs-cont is equiv to abs-cont/ mu_g+/- of its associated measures}
Under the hypotheses of Theorem~\ref{thm: F is ACg is equiv to mu_F satisfying 2 properties}, let us set
\[
\overline{\mu_F}(E):=\mu_F (E\cap A_g^+) \text{ and } \underline{\mu_F}(E):=\mu_F (E\cap A_g^-),
\]
for every $\widetilde g$\--measurable set $E$. Then, the following statements are equivalent:
\begin{enumerate}
\item $F \in \mathcal{AC}_{g}([a,b],{\mathbb R})$;
\item The measures $\overline{\mu_F}$ and $\underline{\mu_F}$ satisfy:
 \begin{enumerate}
 \item $\overline{\mu_F}\ll \mu_g^+$;
 \item $\underline{\mu_F}\ll \mu_g^-$.
 \end{enumerate}
\end{enumerate}
\end{cor}
\begin{rem}
In Corollary~\ref{cor:f g-abs-cont is equiv to abs-cont/ mu_g+/- of its associated measures}, notice that $\mu_F=\overline{\mu_F}+\underline{\mu_F}$, and we have, in general,
\[
\overline{\mu_F} \neq \mu_F^+\text{\quad and \quad} \underline{\mu_F} \neq -\mu_F^-,
\]
where $\mu_F^+$ and $\mu_F^-$ are the measures given by Hahn's decomposition $\mu_F=\mu_F^+ - \mu_F^-$. This remark follows immediately from the Hahn decomposition of the interval $[a,b]$ under $\mu_g$ and $\mu_F$. Indeed, consider, for instance, $[a,b]=[0,2]$, and the function $g:{[0,2]}\to{\mathbb R}$ given by~\eqref{eq:rem:topologies comparison: g}. We have that $A_g^+=[0,1]$ and $A_g^-=(1,2]$. Now, let us consider the function $F:[0,2]\to{\mathbb R}$ defined by $F(t)=-g(t)$ for all $t\in[0,2]$. We have that $F\in \operatorname{\operatorname{BV}^{lc}}([0,2],\mathbb{R})$. Moreover, $A_F^+=(1,2]$ and $A_F^-=[0,1]$. Observe that
\[
\overline{\mu_F}((1,2]):=\mu_F ((1,2]\cap A_g^+)=0\neq 1=F(2)-F(1)=\mu_F^+((1,2]),
\]
and
\[
\underline{\mu_F}([0,1]):=\mu_F ([0,1]\cap A_g^-)=0\neq -1=-(F(0)-F(1))=-\mu_F^-([0,1]).
\]
\end{rem}
In the next lemma, we demonstrate that one can construct $g$\--absolutely continuous functions from $g$\--integrable functions.
\begin{lem}\label{lem:LS-primitive of a func in L1g is g-abs-cont}
 Let $f\in L^1_g([a,b),{\mathbb R})$, and set $F(t):=\int_{[a,t)}f \,\operatorname{d} \mu_g$. Then $F\in \mathcal{AC}_{g}([a,b],{\mathbb R})$.
\end{lem}

\begin{proof}
 We prove this result for $f\geqslant 0$ $|\mu_g|$\--almost everywhere, since the general case is a difference of two non\--negative $\mu_g$\--integrable functions. Let $\varepsilon>0$ be fixed. Since $f\in L^{1}_g[a,b),\mathbb{R})$, there exists $\delta>0$ such that
 \[
 \int_{E} f\, \operatorname{d} |\mu_g|<\varepsilon,
 \]
 for every $\widetilde{g}$\--measurable set $E$ such that $|\mu_g|(E)<\delta$. Let us consider a family $\{(a_i , b_i)\}_{i=1}^{i=m}$ of pairwise disjoint open subintervals of $[a,b]$ such that
 \[
 \sum_{i=1}^{m} \operatorname{var}_g[a_i,b_i] <\delta.
 \]
 Let $E_:=\bigcup_{i=1}^{m}[a_i,b_i)$, then
\[
 |\mu_g|(E)=|\mu_g|\left( \bigcup_{i=1}^{m}[a_i,b_i)\right) = \sum_{i=1}^{m} |\mu_g|([a_i,b_i))=\sum_{i=1}^{m} \operatorname{var}_g[a_i,b_i] <\delta.
\]
 Using the definition of $F$, we obtain
\[
 \sum_{i=1}^{m} |F(b_i)-F(a_i)|\leqslant\sum_{i=1}^{m} \int_{[a_i,b_i)} f\, \operatorname{d} |\mu_g|= \int_E f\, \operatorname{d} |\mu_g|<\varepsilon.
\]
Therefore, $F\in \mathcal{AC}_{g}([a,b],{\mathbb R})$.
\end{proof}

\section{Generalized Fundamental Theorem of Calculus}
We devote this section for generalized versions of the Fundamental Theorem of Calculus involving non\--monotonic derivators. Throughout this section, we continue to consider a closed interval $[a,b]\subset \mathbb{R}$ and a non\--constant derivator $g:[a,b] \to \mathbb{R}$ in $\operatorname{BV}^{\operatorname{lc}}([a,b],{\mathbb R})$. We start with the following result which is the first part of the fundamental theorem of calculus (derivative of the integral), a generalization of Theorem~\ref{derintstieljtes} for nondecreasing derivators and, of \cite[Theorem~6.3]{FrTo} for the case of non\--monotonic derivators with controlled variation.
\begin{thm}\label{lem:f in L1g yields F'g=f}
 Let $g\in\operatorname{BV}^{\operatorname{lc}}([a,b],{\mathbb R})$ be non\--constant, $f\in L^1_g([a,b),\mathbb{R})$, and $F(t):=\int_{[a,t)}f \,\operatorname{d}\mu_g$. Then $F_g'=f$ $|\mu_g|$\--a.e in $[a,b]$.
\end{thm}
\begin{proof}

If $t\in [a,b)\cap D_g$, then we have that
\[
F(t^+)-F(t)=\int_{\{t\}}f\, \operatorname{d} \mu_g=f(t)\mu_g(\{t\})=f(t)(g(t^+)-g(t)).
\]
Thus, using Remark~\ref{rem:ld} we obtain that $F_g'(t)=f(t)$. Now, let us prove that $F_g'=f$ for $|\mu_g|$\--a.e. $t\in [a,b]\setminus D_g$. First, observe that if we set $\mathcal{N}:=N_g\cup C_g=N_{\widetilde{g}}\cup C_{\widetilde{g}}$, then $\mu_{\widetilde{g}}(\mathcal{N})=0$ since $\mu_{\widetilde{g}}(N_{\widetilde{g}})=\mu_{\widetilde{g}}(C_{\widetilde{g}})=0$ according to~\cite[Propositions~2.5 and~2.6]{PR}. Thus, it suffices to prove that $F_g'=f$ for $|\mu_g|$\--a.e. $t\in [a,b]\setminus (\mathcal{N}\cup D_{g})$. In addition, notice that if $t\notin \mathcal{N}$, then $\widetilde{g}(s)\neq \widetilde{g}(t)$ for every $s\neq t$.

Write $g=g_1-g_2+{g(a)}$ as in Corollary~\ref{cor:relating mu_g+ and mu_g1, mu_g- and mu_g2} and recall that $|\mu_g|=\mu_{\widetilde{g}}$. For $|\mu_g|$\--a.e. $t\in [a,b]\setminus (\mathcal{N}\cup D_g)$ and $s$ sufficiently close to $t$, consider the notation
\begin{equation*}\label{eq:notation[|s,t|)}
\llbracket s,t \rrparenthesis:=[\min\{s,t\},\max\{s,t\}).
\end{equation*}
 Thus,
\begin{align*}
 \frac{\int_{[a,s)} \chi_{A_g^+}\,\operatorname{d} \mu_{\widetilde{g}} -\int_{[a,t)} \chi_{A_g^+}\,\operatorname{d} \mu_{\widetilde{g}}}{\widetilde{g}(s)-\widetilde{g}(t)}& =\operatorname{sgn}(s-t)\frac{\int_{\llbracket s,t \rrparenthesis} \chi_{A_g^+}\,\operatorname{d} \mu_{\widetilde{g}}}{\widetilde{g}(s)-\widetilde{g}(t)}=\operatorname{sgn}(s-t)\frac{|\mu_g|\left(\llbracket s,t \rrparenthesis \cap A_g^+\right)}{\widetilde{g}(s)-\widetilde{g}(t)}\\
 &=\operatorname{sgn}(s-t)\frac{\mu_g^+\big(\llbracket s,t \rrparenthesis\big) }{\widetilde{g}(s)-\widetilde{g}(t)}
 =\frac{g_1(s)-g_1(t)}{\widetilde{g}(s)-\widetilde{g}(t)}.
\end{align*}
Since $\widetilde{g}$ is left\--continuous and nondecreasing, {it follows from Theorem~\ref{derintstieljtes}, that for $\mu_{\widetilde{g}}$\--a.e. $t\in [a,b]\setminus (\mathcal{N}\cup D_g)$ }
\begin{align*}
\chi_{A_g^+}(t)&=\frac{\operatorname{d}}{\operatorname{d}_{\widetilde{g}}t} \int_{[a,t)} \chi_{A_g^+}\,\operatorname{d} \mu_{\widetilde{g}}=\lim\limits_{s\to t}\frac{\int_{[a,s)} \chi_{A_g^+}\,\operatorname{d} \mu_{\widetilde{g}} -\int_{[a,t)} \chi_{A_g^+}\,\operatorname{d} \mu_{\widetilde{g}}}{\widetilde{g}(s)-\widetilde{g}(t)}=\lim\limits_{s\to t}\frac{g_1(s)-g_1(t)}{\widetilde{g}(s)-\widetilde{g}(t)}.
\end{align*}
This implies that
\[
\lim\limits_{s\to t}\frac{g_1(s)-g_1(t)}{\widetilde{g}(s)-\widetilde{g}(t)}=
\begin{dcases}
1, & \mbox{for $\mu_{\widetilde{g}}$\--a.e. }t\in A_g^+\setminus (\mathcal{N}\cup D_g),\\
0, & \mbox{for $\mu_{\widetilde{g}}$\--a.e. }t\in A_g^-\setminus (\mathcal{N}\cup D_g).
\end{dcases}
\]
Similarly, we obtain that
\[
\lim\limits_{s\to t}\frac{g_2(s)-g_2(t)}{\widetilde{g}(s)-\widetilde{g}(t)}=
\begin{dcases}
1, & \mbox{for $\mu_{\widetilde{g}}$\--a.e. }t\in A_g^-\setminus (\mathcal{N}\cup D_g),\\
0, & \mbox{for $\mu_{\widetilde{g}}$\--a.e. }t\in A_g^+\setminus (\mathcal{N}\cup D_g).
\end{dcases}
\]
Since $g(s)-g(t)=g_1(s)-g_1(t)-(g_2(s)-g_2(t))$ for all $t,s\in[a,b]$, we obtain
\[
\lim_{\substack{s\to t\\g(s)\neq g(t)}}\frac{g(s)-g(t)}{\widetilde{g}(s)-\widetilde{g}(t)}=
\begin{dcases}
	\lim\limits_{s\to t}\frac{g_1(s)-g_1(t)}{\widetilde{g}(s)-\widetilde{g}(t)}=1 & \mbox{for $|\mu_g|$\--a.e. }t\in A_g^+\setminus (\mathcal{N}\cup D_g),\\
	-\lim\limits_{s\to t}\frac{g_2(s)-g_2(t)}{\widetilde{g}(s)-\widetilde{g}(t)}=-1 & \mbox{for $|\mu_g|$\--a.e. } t\in A_g^-\setminus (\mathcal{N}\cup D_g).
\end{dcases}
\]
Given that $f\in L^1_g([a,b),{\mathbb R})$, Lemmata~\ref{lem:g-abs-cont func has bounded var} and~\ref{lem:LS-primitive of a func in L1g is g-abs-cont} imply that $F\in \mathcal{AC}_g([a,b],\mathbb{R})\subset\operatorname{\operatorname{BV}^{lc}}([a,b],{\mathbb R})$.
 We define,
\[
F_1(t):=\int_{[a,t)}f\, \operatorname{d} \mu_{g_1},\quad\text{and}\quad F_2(t):=\int_{[a,t)}f\, \operatorname{d} \mu_{g_2}\quad \text{ for every $t\in [a,b]$.}
\]
By construction, the functions $F_1,F_2:[a,b]\to \mathbb{R}$ are left-continuous, and using $\mu_g=\mu_{g_1}-\mu_{g_2}$ we obtain for all $t\in[a,b]$,
\[
F(t)=\int_{[a,t)}f\, \operatorname{d} \mu_g=F_1(t)-F_2(t).
\]
Hence, by Theorem~\ref{derintstieljtes}, the $g_k$\--derivative of $F_k$ exists $\mu_{g_k}$\--a.e. for $k=1,2$
\[
(F_1)_{g_1}'=f \,\mu_{g_1}\text{\--a.e. in $[a,b)$} \quad\text{and}\quad (F_2)_{g_2}'=f \, \mu_{g_2}\text{\--a.e. in $[a,b)$.}
\]
By Corollary~\ref{cor:relating mu_g+ and mu_g1, mu_g- and mu_g2}, $\mu_{g_1}=\mu_g^+$, $\mu_{g_2}=\mu_g^-$, and $|\mu_g|= \mu_g^+ + \mu_g^-$, so it follows that
\[
(F_1)_{g_1}'=f \,\text{ $|\mu_g|$\--a.e. in $A_g^+$} \quad\text{and}\quad (F_2)_{g_2}'=f \, \text{$|\mu_g|$\--a.e. in $A_g^-$.}
\]
Furthermore, for $|\mu_g|$\--almost every $t\in [a,b]\setminus (\mathcal{N}\cup D_g)$,
\begin{align*}
 \lim\limits_{s\to t} \frac{F_1(s)-F_1(t)}{\widetilde{g}(s)-\widetilde{g}(t)} = & \lim\limits_{s\to t} \frac{\int_{[a,s)}f\, \operatorname{d} \mu_{g_1}- \int_{[a,t)}f\, \operatorname{d} \mu_{g_1} }{\widetilde{g}(s)-\widetilde{g}(t)}=\lim\limits_{s\to t} \frac{\int_{[a,s)}f\chi_{A_g^+}\, \operatorname{d} \mu_{\widetilde{g}}- \int_{[a,t)}f\chi_{A_g^+}\, \operatorname{d} \mu_{\widetilde{g}} }{\widetilde{g}(s)-\widetilde{g}(t)} \\
 = & f(t)\chi_{A_g^+}(t),
\end{align*}
and, similarly,
\[
 \lim\limits_{s\to t} \frac{F_2(s)-F_2(t)}{\widetilde{g}(s)-\widetilde{g}(t)}= f(t)\chi_{A_g^-}(t).
\]

For a $t\in A_g^+ \setminus (\mathcal{N}\cup D_g)$ such that
\[  \lim\limits_{s\to t}\frac{g_1(s)-g_1(t)}{\widetilde{g}(s)-\widetilde{g}(t)}=\lim\limits_{\substack{s\to t \\g(s) \neq g(t)}}\frac{g(s)-g(t)}{\widetilde{g}(s)-\widetilde{g}(t)}=1,\]
(recall this condition is satisfied for $|\mu_g|$\--almost every $t\in A_g^+ \setminus (\mathcal{N}\cup D_g)$) let us consider the function $h:[a,b]\to{\mathbb R}$ defined by
\[
h(s)=\begin{dcases}
 \frac{F(s)-F(t^*)}{g(s)-g(t^*)}-f(t), & \mbox{if } g(s)\neq g(t^*), \\
 0, & \text{otherwise.}
 \end{dcases}
\]
Observe that, for $|\mu_g|$\--almost every $t\in A_g^+ \setminus (\mathcal{N}\cup D_g)$, $t^*=t$, $F(t)=F(t^*)$ and $g(t)=g(t^*)$. Then,
\begin{align*}
\lim\limits_{\substack{s\to t \\g(s) \neq g(t)}} h(s) =& \lim\limits_{\substack{s\to t \\g(s) \neq g(t)}} \left( \frac{F(s)-F(t)}{g(s)-g(t)}-f(t)\right)= \lim\limits_{\substack{s\to t \\g(s) \neq g(t)}} \left(\frac{F_1(s)-F_1(t)}{g(s)-g(t)}-\frac{F_2(s)-F_2(t)}{g(s)-g(t)}-f(t)\right)\\
 =& \lim\limits_{\substack{s\to t \\{g(s) \neq g(t)}}} \frac{F_1(s)-F_1(t)}{g_1(s)-g_1(t)}\frac{g_1(s)-g_1(t)}{\widetilde{g}(s)-\widetilde{g}(t)}\frac{\widetilde{g}(s)-\widetilde{g}(t)}{g(s)-g(t)} -\frac{F_2(s)-F_2(t)}{\widetilde{g}(s)-\widetilde{g}(t)}\frac{\widetilde{g}(s)-\widetilde{g}(t)}{g(s)-g(t)} -f(t)\\ =&(F_1)'_{g_1}(t)-f(t)=0=h(t),
\end{align*}
where the limits involved are taken for the points $s\in[a,b]$ such that $g(s)\ne g(t^*)$ (for the case $g(s)=g(t^*)$, it is trivial, as $h(s)=0$).
Thus, $h$ is continuous at such $t$. Hence, $F$ is $g$\--differentiable at $t$. Arguing similarly for $t\in A_g^-\setminus (\mathcal{N}\cup D_g)$ such that
\[  \lim\limits_{s\to t}\frac{g_2(s)-g_2(t)}{\widetilde{g}(s)-\widetilde{g}(t)}=-\lim\limits_{\substack{s\to t \\g(s) \neq g(t)}}\frac{g(s)-g(t)}{\widetilde{g}(s)-\widetilde{g}(t)}=1,\]
we deduce that $h$ is continuous at $t$. In conclusion:
\begin{itemize}
 \item $F_g'$ exists $|\mu_g|$\--a.e. in $A_g^+$ and $F_g'=(F_1)_{g_1}'=f$.
 \item $F_g'$ exists $|\mu_g|$\--a.e. in $A_g^-$ and $F_g'=(F_2)_{g_2}'=f$.
\end{itemize}
Hence, $F_g'=f$ $|\mu_g|$\--a.e. in $[a,b]$.
\end{proof}
Now, we state the Fundamental Theorem of Calculus for Lebesgue\--Stieltjes integrals (for the integral of the Stieltjes derivative).
\begin{thm}[Fundamental Theorem of Calculus for Lebesgue\--Stieltjes integrals]\label{thm:FTC1}
Let $g\in\operatorname{BV}^{\operatorname{lc}}([a,b],{\mathbb R})$ be non\--constant. Then, the following statements are equivalent:
\begin{enumerate}
 \item $F \in \mathcal{AC}_{g}([a,b],{\mathbb R})$;
 \item \begin{enumerate}
 \item $F_g'$ exists $|\mu_g|$\--a.e. in $[a,b]$;
 \item $F_g' \in L^1_g([a,b),{\mathbb R})$;
 \item $F(t)=F(a)+ \int_{[a,t)}F_g'\, \operatorname{d} \mu_g,$ for $|\mu_g|$\--a.e. $t\in [a,b]$.
 \end{enumerate}
\end{enumerate}
\end{thm}

\begin{proof}
$(2)\Rightarrow(1)$. This implication follows immediately from Lemma~\ref{lem:LS-primitive of a func in L1g is g-abs-cont}.

$(1)\Rightarrow(2)$. If $F\in \mathcal{AC}_{g}([a,b],{\mathbb R})$, then, using Corollary~\ref{cor:f g-abs-cont is equiv to abs-cont/ mu_g+/- of its associated measures}, we obtain
 \[  \overline{\mu_F}\ll \mu_g^+ \text{ and } \underline{\mu_F}\ll \mu_g^-.
 \]
By applying the Radon\--Nikod\'{y}m Theorem~\cite[Theorem~5.3.2]{BenedettoCzaja2009integration}, there exist unique functions $f^+\in L^{1}_{\mu_g^+}([a,b),\mathbb{R})$ and $f^-\in L^{1}_{\mu_g^-}([a,b),\mathbb{R})$ such that
\[
 \overline{\mu_F}(E^+)=\int_{E^+}f^+ \, \operatorname{d} \mu_g^+,
\]
and
\[
\underline{\mu_F}(E^-)=\int_{E^-}f^- \, \operatorname{d} \mu_g^-,
\]
for any $\mu_g^+$\--measurable set $E^+ \subset A_g^+$ and $\mu_g^-$\--measurable set $E^- \subset A_g^-$.

Let us define the function $f:[a,b] \to {\mathbb R}$ by
\[
f(t)=\begin{dcases}
 f^+(t), & \mbox{if } t\in A_g^+,\\
 -f^-(t), & \mbox{if } t\in A_g^-.
 \end{dcases}
\]
It follows from Theorem~\ref{thm:consequence:Hahn+Jordan Decomp-sets A+g & A-g} that $f\in L^1_g([a,b),{\mathbb R})$, since
\begin{align*}
\int_{[a,b)}|f| \, \operatorname{d} |\mu_g| =\int_{[a,b)\cap A_g^+}|f| \, \operatorname{d} |\mu_g| + \int_{[a,b)\cap A_g^-}|f| \, \operatorname{d} |\mu_g|
 = \int_{[a,b)}|f^+| \, \operatorname{d} \mu_g^+ + \int_{[a,b)}|f^-| \, \operatorname{d} \mu_g^-
 <\infty.
\end{align*}
Moreover, for a $\widetilde{g}$\--measurable set $E \subset [a,b)$, we have that
\begin{align*}
 \mu_F(E)& =\mu_F(E\cap A_g^+)+ \mu_F(E\cap A_g^-)=\overline{\mu_F}(E\cap A_g^+)+\underline{\mu_F}(E\cap A_g^-) =\int_{E\cap A_g^+} f^+ \, \operatorname{d} \mu_g^+ + \int_{E\cap A_g^-} f^- \, \operatorname{d} \mu_g^-\\
 & =\int_{E\cap A_g^+} f^+ \, \operatorname{d} \mu_g - \int_{E\cap A_g^-} f^- \, \operatorname{d} \mu_g= \int_{E} f \, \operatorname{d} \mu_g.
\end{align*}
In particular, for $E=[a,t)$, with $t\in [a,b]$, we obtain
\[
F(t)-F(a)=\mu_F([a,t))=\int_{[a,t)} f \operatorname{d} \mu_g.
\]
Hence, using Theorem~\ref{lem:f in L1g yields F'g=f}, we conclude that $F_g'=f$ $|\mu_g|$\--a.e. in $[a,b]$.
\end{proof}
The following result strengthens the conclusion of Theorem~\ref{lem:f in L1g yields F'g=f} under stronger regularity conditions, even in the context of monotonic derivators. It allows the computation of the Stieltjes derivative at each point, based on the revised definition of the Stieltjes derivative through the notation $t^*$ in~\eqref{eq:notation of t^*}. Additionally, classical literature, such as in~\cite[Proposition~A.2.8(iv)]{AL}, asserts that a primitive~$F$ of a continuous function~$f$ is always differentiable, and $F'=f$ holds everywhere. In the context of derivators of bounded variation, interestingly, we find that the Stieltjes analogous version necessitates the introduction of an additional assumption\---one that is always satisfied in the monotonic case and, in particular, holds in the classical setting.
\begin{thm}\label{thmfcf} Let $g\in \operatorname{BV}^{\operatorname{lc}}([a,b],{\mathbb R})$ be such that $b\notin N_g^+$.
 Let us define the function $\varphi:[a,b]\to{\mathbb R}$ by
\begin{equation}\label{eq:thm:varphi function}
\varphi(t):=\begin{dcases}
 \liminf_{s\to t}\left|\frac{g(s)-g(t)}{\widetilde{g}(s)-\widetilde{g}(t)}\right| & \mbox{if } t^*\in[a,b]\setminus(D_g \cup N_g), \\
 \liminf_{s\to t^{*+}}\left|\frac{g(s)-g(t^*)}{\widetilde{g}(s)-\widetilde{g}(t^*)}\right| & \mbox{if } t^*\in(D_g\cup N_g^+),\\
 \liminf_{s\to t^-}\left|\frac{g(s)-g(t)}{\widetilde{g}(s)-\widetilde{g}(t)}\right| & \mbox{if } t^*\in N_g^- ,
 \end{dcases}
\end{equation}
and assume that
\begin{equation}\label{eq:thm:varphi >0}
 \varphi(t)>0 \quad \text{for all $t\in[a,b]$.}
\end{equation}
 Let $f\in {\mathcal C}_g([a,b],{\mathbb R})\cap L^1_g([a,b),\mathbb{R})$ and $F(t):=\int_{[a,t)} f\, \operatorname{d} \mu_g$ for $t\in[a,b]$. Then the following statements hold:
\begin{enumerate}
 \item $F \in \mathcal{AC}_{g}([a,b],{\mathbb R})$.
 \item $F_g'(t)=f(t)$ for all $t\in [a,b]$.
\end{enumerate}
\end{thm}
\begin{proof}
 1. Since $f \in L^1_{g}([a,b),\mathbb{R})$, then, according to Lemma~\ref{lem:LS-primitive of a func in L1g is g-abs-cont}, $F \in \mathcal{AC}_{g}([a,b],{\mathbb R})$.

2. Fix $t\in [a,b]$ and $\varepsilon>0$. We distinguish four cases:

 \textbf{Case 1:} $t^* \in[a,b]\setminus (D_g\cup N_g)$, then $t^*=t$ and $t\in[a,b]\setminus(D_g\cup C_g \cup N_g)$.
In this case, since \[
\varphi(t)=\liminf_{s\to t}\left|\frac{g(s)-g(t)}{\widetilde{g}(s)-\widetilde{g}(t)}\right|>0,
\]
 there exists $\delta_1>0$ such that, for $s\in [a,b]$ satisfying $|s-t|<\delta_1$, we have
\[
\frac{\varphi(t)}{2}< \left|\frac{g(s)-g(t)}{\widetilde{g}(s)-\widetilde{g}(t)}\right|,
\]
which implies that
\begin{equation}\label{eq:proof:case1:comparing g and tilde(g) around a point}
|\widetilde{g}(s)-\widetilde{g}(t)|<\frac{2}{\varphi(t)} |g(s)-g(t)|.
\end{equation}
Since $f$ is $g$\--continuous, there exists $\delta_2\in{\mathbb R}^+$ such that
\[
|f(s)-f(t)|<\xi:=\frac{\varepsilon\varphi(t)}{2}, \text{ for } s \in [a,b] \text{ such that } |\widetilde g(s)-\widetilde g(t)|<\delta_2.
\]
Given that $t\not\in D_g$, $\widetilde g$ is continuous at $t$, so there exists $\delta\in(0,\delta_1]$ such that, if $|s-t|<\delta$, then $|\widetilde g(s)-\widetilde g(t)|<\delta_2$. Thus, if $s\in(t-\delta,t+\delta)\cap[a,b]$, we obtain
\[
f(t)-\xi< f(s) < f(t)+\xi.
\]
Since $\mu_g=\mu_g^+ -\mu_g^-$, {then} for all $s\in(t-\delta,t+\delta)\cap[a,b]$ such that $s>t$ we have that
\begin{equation}\label{eqd1}
(f(t)-\xi)\mu_g^+([t,s))\leqslant \int_{[t,s)} f\, \operatorname{d} \mu_g^+ \leqslant (f(t)+\xi)\mu_g^+([t,s)),
\end{equation}
and
\[
(f(t)-\xi)\mu_g^-([t,s))\leqslant \int_{[t,s)} f\, \operatorname{d} \mu_g^- \leqslant (f(t)+\xi)\mu_g^-([t,s)),
\]
or equivalently,
\begin{equation}\label{eqd2}
-(f(t)+\xi)\mu_g^-([t,s))\leqslant -\int_{[t,s)} f\, \operatorname{d} \mu_g^- \leqslant -(f(t)-\xi)\mu_g^-([t,s)).
\end{equation}
Given that
\[
F(s)-F(t)=\int_{[t,s)} f\, \operatorname{d} \mu_g=\int_{[t,s)} f\, \operatorname{d} \mu_g^+-\int_{[t,s)} f\, \operatorname{d} \mu_g^-,
\]
adding~\eqref{eqd1} and~\eqref{eqd2}, we get
\[  (f(t)-\xi)\mu_g^+([t,s))-(f(t)+\xi)\mu_g^-([t,s))\leqslant F(s)-F(t) \leqslant (f(t)+\xi)\mu_g^+([t,s)) -(f(t)-\xi)\mu_g^-([t,s)),\]
that is,
\[  \mu_g([t,s))f(t)-\xi|\mu_g|([t,s))\leqslant F(s)-F(t) \leqslant \mu_g([t,s))f(t)+\xi|\mu_g|([t,s)).\]
Hence,
\[  -\xi|\mu_g|([t,s))\leqslant F(s)-F(t)-\mu_g([t,s))f(t) \leqslant\xi|\mu_g|([t,s)),\]
or, equivalently,
\begin{equation}\label{eq:thm Prelim:ineq for t>t0}
	|F(s)-F(t)-f(t)(g(s)-g(t))|\leqslant\xi|\mu_g|([t,s))=\xi(\widetilde{g}(s)-\widetilde{g}(t)).\end{equation}
Similarly, for all $s\in(t-\delta,t+\delta)\cap[a,b]$ such that $s<t$, we obtain
\begin{equation*}
\mu_g([s,t))f(t)-\xi|\mu_g|([s,t))\leqslant F(t)-F(s) \leqslant \mu_g([s,t))f(t)+\xi|\mu_g|([s,t)).
\end{equation*}
Therefore,
\[
-\mu_g([s,t))f(t)-\xi|\mu_g|([s,t))\leqslant F(s)-F(t) \leqslant -\mu_g([s,t))f(t)+\xi|\mu_g|([s,t)),\]
and, thus,
\begin{equation}\label{eq:thm Prelim:ineq for t<t0}
	|F(s)-F(t)-f(t)(g(s)-g(t))|\leqslant\xi|\mu_g|([s,t))=\xi(\widetilde{g}(t)-\widetilde{g}(s)).\end{equation}
Given that $\widetilde g$ is nondecreasing, using~\eqref{eq:thm Prelim:ineq for t>t0} and~\eqref{eq:thm Prelim:ineq for t<t0}, we deduce that, for all $s\in(t-\delta,t+\delta)\cap[a,b]$,
\[
|F(s)-F(t)-f(t)(g(s)-g(t))|\leqslant \xi|\widetilde{g}(s)-\widetilde{g}(t)|.
\]
Combining this inequality with~\eqref{eq:proof:case1:comparing g and tilde(g) around a point} we obtain, for all $s\in(t-\delta,t+\delta)\cap[a,b]$,
\begin{equation}\label{eq:proof:Case1:|F(s)-F(t)-f(t)(g(s)-g(t))|<e|g(s)-g(t)|}
| F(s)-F(t)-f(t)(g(s)-g(t))| \leqslant \xi|\widetilde{g}(s)-\widetilde{g}(t)|<\varepsilon |g(s)-g(t)|.
\end{equation}
\textbf{Case 2:} $t^*\in N_g^+$. In this case, since
\[
\varphi(t)=\liminf_{s\to t^{*+}}\left|\frac{g(s)-g(t^*)}{\widetilde{g}(s)-\widetilde{g}(t^*)}\right|>0,
\]
there exists $\delta_1>0$ such that, for $s\in [a,b]$ satisfying $|s-t^*|<\delta_1$, we have
\[
\frac{\varphi(t)}{2}< \left|\frac{g(s)-g(t^*)}{\widetilde{g}(s)-\widetilde{g}(t^*)}\right|,
\]
which implies that
\begin{equation}\label{eq:proof:case2:comparing g and tilde(g) around a point}
	|\widetilde{g}(s)-\widetilde{g}(t^*)|<\frac{2}{\varphi(t)} |g(s)-g(t^*)|.
\end{equation}
Now, by Proposition~\ref{prop:g-cont fonc properties}, we have that $f(t^*)=f(t)$. Since $f$ is $g$\--continuous at $t^*$, there exists $\delta_2>0$ such that
\[
|f(s)-f(t)|<\xi:=\frac{\varepsilon\varphi(t)}{2}, \text{ for } s \in [a,b] \text{ such that } |\widetilde g(s)-\widetilde g(t^*)|<\delta_2.
\]
Given that $t^*\in N_g^+$, $\widetilde g$ is continuous at $t^*$, so there exists $\delta\in(0,\delta_1]$ such that, if $|s-t^*|<\delta$, then $|\widetilde g(s)-\widetilde g(t^*)|<\delta_2$. Thus, if $s\in(t^*,t^* +\delta)\cap[a,b]$, we obtain,
\[
f(t)-\xi< f(s) < f(t)+\xi.
\]
Arguing similarly to the previous case, we deduce that, for all $s\in(t^*,t^* +\delta)\cap[a,b]$,
\[
|F(s)-F(t^*)-f(t)(g(s)-g(t^*))|\leqslant\xi|\mu_g|([t^*,s))=\xi(\widetilde{g}(s)-\widetilde{g}(t^*)).
\]
Using~\eqref{eq:proof:case2:comparing g and tilde(g) around a point}, we get
\begin{equation}\label{eq:proof:Case2:|F(s)-F(t)-f(t)(g(s)-g(t))|<e|g(s)-g(t)|}
| F(s)-F(t^*)-f(t)(g(s)-g(t^*))| \leqslant \xi|\widetilde{g}(s)-\widetilde{g}(t^*)|<\varepsilon|g(s)-g(t^*)|.
\end{equation}
\textbf{Case 3:} $t^*\in N_g^-$, then $t^*=t$. In this case, since
\[
\varphi(t)=\liminf_{s\to t^{-}}\left|\frac{g(s)-g(t)}{\widetilde{g}(s)-\widetilde{g}(t)}\right|>0,
\]
there exists $\delta_1>0$ such that, for $s\in [a,b]\cap (t-\delta_1,t)$ , we have
\[
\frac{\varphi(t)}{2}< \left|\frac{g(s)-g(t)}{\widetilde{g}(s)-\widetilde{g}(t)}\right|,
\]
which implies that
\begin{equation}\label{eq:proof:case3:comparing g and tilde(g) around a point}
	|\widetilde{g}(s)-\widetilde{g}(t)|<\frac{2}{\varphi(t)} |g(s)-g(t)|.
\end{equation}
Arguing analogously to {\bf Case~2} and using~\eqref{eq:proof:case3:comparing g and tilde(g) around a point}, we obtain that there exists $\delta\in{\mathbb R}^+$ such that, if $s\in [a,b]\cap (t-\delta,t)$, then
\begin{equation}\label{eq:proof:Case3:|F(s)-F(t)-f(t)(g(s)-g(t))|<e|g(s)-g(t)|}
| F(s)-F(t)-f(t)(g(s)-g(t))| < \varepsilon |g(s)-g(t)|.
\end{equation}
\textbf{Case 4:} $t^*\in D_g$. In this case, $g$ is discontinuous at~$t^*$, and we have either $g(t^{*+})>g(t^*)$ or $g(t^{*+})<g(t^*)$. Since the variation function~$\widetilde{g}$ of~$g$ is nondecreasing, it follows that $\widetilde{g}(t^{*+})>\widetilde{g}(t^*)$, and by Remark~\ref{rem:Dg+(t)=Dtil(g)+(t) + var[a,b]=til(g)(b)-tilde(a)} we obtain
\[
|g(t^{*+})-g(t^*)|=\widetilde{g}(t^{*+})-\widetilde{g}(t^*).
\]
Thus, we have that
\[
\varphi(t)=\liminf_{s\to t^{*+}}\left|\frac{g(s)-g(t^*)}{\widetilde{g}(s)-\widetilde{g}(t^*)}\right|=\frac{|\Delta^+ g(t^*)|}{\Delta^+ \widetilde g(t^*)}=1>0.
\]
Consequently, there exists $\delta_1 \in\mathbb{R}^+$ such that for $s\in (t^*,t^*+\delta_1)$, we have
\[
\frac{1}{2}<\left|\frac{g(s)-g(t^*)}{\widetilde{g}(s)-\widetilde{g}(t^*)}\right|,
\]
which implies that
\begin{equation}\label{eq:proof:case4:comparing g and tilde(g) around a point}
	|\widetilde{g}(s)-\widetilde{g}(t^*)|<2|g(s)-g(t^*)|.
\end{equation}
Let $M:=\max\left\{1,|f(t)|\right\}$. Since $f$ is $g$\--integrable on $[a,b)$, there exists $\eta>0$ so that for $s\in(t^*,b]$ such that $|\widetilde{g}(s)-\widetilde{g}(t^{*+})|=|\mu_g|((t^*,s))<\eta$, we have
\[
 \int_{(t^*,s)}|f|\operatorname{d}|\mu_g|< \frac{\varepsilon}{4}\Delta^+\widetilde g(t^*).
\]
Additionally, since there exists $\widetilde g(t^{*+})$, there is $\delta\in(0,\delta_1]$ such that for all $s\in(t^*,t^* +\delta)\cap [a,b]$, we have
\[  |\widetilde g(s)-\widetilde g(t^{*+})|<\min\left\{\frac{\varepsilon}{4M}\Delta^+\widetilde g(t^*),\eta\right\}.\]
By Proposition~\ref{prop:g-cont fonc properties}, we have that $f(t)=f(t^*)$. So, for $s\in(t^*,t^{*}+\delta)\cap [a,b]$,
\begin{align*}
\left|F(s)-F(t^*)-f(t)(g(s)-g(t^*))\right|=& \left|\int_{[t^*,s)}f\operatorname{d}\mu_g-f(t^*)(g(s)-g(t^*))\right|\\
=&\left|\int_{\{t^*\}}f\operatorname{d}\mu_{g}+\int_{(t^*,s)}f\operatorname{d}\mu_{g}-f(t^*)(g(s)-g(t^*))\right|\\
=& \left|f(t^*)(g(t^{*+})-g(t^*))+\int_{(t^*,s)}f\operatorname{d}\mu_g-f(t^*)(g(s)-g(t^*))\right|\\
=& \left|f(t^*)(g(t^{*+})-g(s))+\int_{(t^*,s)}f\operatorname{d}\mu_g\right|\\
\leqslant &|f(t)||g(s)-g(t^{*+})|+ \int_{(t^*,s)}|f|\operatorname{d}|\mu_g|\\
< &|f(t)||\widetilde{g}(s)-\widetilde{g}(t^{*+})|+ \frac{\varepsilon}{4}\Delta^+\widetilde g(t^*)\\
\leqslant &|f(t)|\frac{\varepsilon}{4M}\Delta^+\widetilde g(t^*)+\frac{\varepsilon}{4}\Delta^+\widetilde g(t^*)\\
\leqslant & \frac{\varepsilon}{2}\Delta^+\widetilde g(t^*)= \frac{\varepsilon}{2}|\widetilde{g}(t^{*+})-\widetilde{g}(t^*))|\\
\leqslant & \frac{\varepsilon}{2}|\widetilde{g}(s)-\widetilde{g}(t^*))|< \varepsilon |g(s)-g(t^*)|.
\end{align*}
The last inequality holds from~\eqref{eq:proof:case4:comparing g and tilde(g) around a point}.

Now let us define the function $h:[a,b]\to\mathbb{R}$ by
\[
h(s):=\begin{dcases}
 \frac{F(s)-F(t^*)}{g(s)-g(t^*)}-f(t), & \mbox{if } g(s)\neq g(t^*), \\
 0, & \mbox{otherwise}.
 \end{dcases}
\]
From~\eqref{eq:proof:Case1:|F(s)-F(t)-f(t)(g(s)-g(t))|<e|g(s)-g(t)|},~\eqref{eq:proof:Case2:|F(s)-F(t)-f(t)(g(s)-g(t))|<e|g(s)-g(t)|}, and~\eqref{eq:proof:Case3:|F(s)-F(t)-f(t)(g(s)-g(t))|<e|g(s)-g(t)|}, it results that there exists $\delta>0$ such that
\[
\begin{aligned}
&|h(s)|<\varepsilon, && \text{for }0<|s-t^*|<\delta,\ \text{with }t^*\not\in D_g\cup N_g,\ \text{and }g(t^*)\ne g(s),\\
&|h(s)|<\varepsilon, && \text{for }0<s-t^*<\delta,\ \text{with }t^*\in D_g\cup N_g^+,\ \text{and } g(t^*)\ne g(s),\\
&|h(s)|<\varepsilon, && \text{for }0<t^*-s<\delta,\ \text{with }t^* \in N_g^-,\ \text{and } g(t^*)\ne g(s).
\end{aligned}
\]
Therefore, $h$ fulfills the assumptions of Definition~\ref{dfn:g-derivative(non monotonic case)}. Hence, $F_g'(t)=f(t)$ for all $t\in [a,b]$.
\end{proof}
\begin{rem}\label{rem:phi(t)>0 is required on continuity points}
In {\bf Case~4} of the proof of Theorem~\ref{thmfcf}, we can clearly observe that the condition $\varphi(t)>0$ holds for all $t\in[a,b]$ such that $t^*\in D_g$. In particular,
 \[
\varphi(t)=\liminf_{s\to t^{+}}\left|\frac{g(s)-g(t)}{\widetilde{g}(s)-\widetilde{g}(t)}\right|=\frac{|\Delta^+ g(t)|}{\Delta^+ \widetilde g(t)}=1>0,\text{ for all }t\in D_g.
\]
Furthermore, if $g$ is nondecreasing, we have that $g=\widetilde{g}+g(a)$, which implies that $\varphi(t)=1$ and hence $\varphi(t)>0$ is always satisfied for all $t\in[a,b]$.
\end{rem}

\begin{rem}It is worthwhile to mention that Theorem~\ref{thmfcf} can be restated with weaker assumptions. In particular, one can notice that in the proof of Theorem~\ref{thmfcf}, the condition for $g$\--continuity of the function $f$ can be limited to the set $[a,b]\setminus (C_g\cup D_g \cup N_g)$. Indeed, observe that $g$\--continuity of $f$ is not necessary on the points of $D_g$. Moreover, it can be observed that on the set $N_g^-$ (resp. $N_g^+$), the requirements can be further restricted to {\it left $g$\--continuity} (resp. {\it right $g$\--continuity}) of $f$ as will be defined in the following point.

In the monotonic case where $g$ is left\--continuous and nondecreasing, the authors in~\cite[Definition~4.4]{FMT2025kernel-Stieltjes-deriv-space} introduced a new space of functions that are $g$\--continuous on $[a,b]\setminus (C_g\cup D_g \cup N_g)$, and left $g$\--continuous from the right (resp. from the left) at the point of $N_g^-$ (resp. $N_g^+$). Analogously to the monotonic case, given $g\in \operatorname{BV}^{\operatorname{lc}}([a,b],{\mathbb R})$, we can define an analogous space which would weaken the $g$\--continuity assumption on $f$ on the set $D_g \cup N_g$ in light of the previous point. Let $f:[a,b]\to\mathbb{R}$ be a function. We say that $f$ is {\it left $g$\--continuous at $t\in (a,b]$} if, for every $\varepsilon>0$, there is $\delta>0$ such that
\[
|f(t)-f(s)|<\varepsilon \text{ for $s\leqslant t$ such that } \widetilde{g}(t)-\widetilde{g}(s)<\delta.
\]
Similarly, we say that $f$ is {\it right $g$\--continuous at $t\in [a,b)$} if, for every $\varepsilon>0$, there is $\delta>0$ such that
\[
|f(t)-f(s)|<\varepsilon \text{ for $s\geqslant t$ such that } \widetilde{g}(s)-\widetilde{g}(t)<\delta.
\]
Now, we can define the space $\mathcal{BD}_g([a,b],\mathbb{R})$ of functions that are bounded, $g$\--continuous on the set $[a,b]\setminus (C_g\cup D_g \cup N_g)$, and left $g$\--continuous from the right (resp. from the left) at the point of $N_g^-$ (resp. $N_g^+$). Observe that, taking int account the previous point, we can weaken the $g$\--continuity assumption on $f$ in Theorem~\ref{thmfcf} by taking $f\in\mathcal{BD}_g([a,b],\mathbb{R})$.
\end{rem}

In the next example, we show that, at least for some derivators, the condition $\varphi(t)>0$ for all $t\in[a,b]$ occurring in Theorem~\ref{thmfcf} is necessary for the assertion $F_g'(t)=f(t)$ to hold for all $t\in [a,b]$.
\begin{exa}\label{exa1} Let $(\alpha_n)_{n\in{\mathbb N}}\subset (0,1)$ be a sequence such that
\begin{equation}\label{eq:sum for psi(x_1)}
\sum_{k=1}^\infty\frac{\alpha_{k+1}}{1+\alpha_{k+1}}\prod_{j=1}^k\frac{1-\alpha_j}{1+\alpha_j}=\frac{1}{2}\frac{ 1-\alpha_1}{1+\alpha_1}.
\end{equation}

In Example~\ref{exa2} we will show that such a choice of the sequence $(\alpha_n)_{n\in{\mathbb N}}$ is possible.	We now define the sequence $(x_n)_{n\in{\mathbb N}}\subset{\mathbb R}$ as follows:
	\[  x_1=1,\quad x_{2n}= \frac{1}{1+\alpha_n}x_{2n-1},\quad x_{2n+1}=(1-\alpha_n)x_{2n},\quad {n\in{\mathbb N}}.\]
	Observe that, $0<x_{n+1}<x_n\leqslant 1$ for every $n\in{\mathbb N}$.	Thus, $x_n\to x_0\in[0,1]$.
	Let us define the function $h:[x_0,1]\to{\mathbb R}$ such that
\[  h(t):=\begin{dcases}-1, & t\in(x_{2n},x_{2n-1}],\ n\in{\mathbb N},\\
		1, & t\in(x_{2n+1},x_{2n}],\ n\in{\mathbb N},\\ 0, & \text{otherwise.}
	\end{dcases}\]
	The function $h$ is Lebesgue integrable. Now let us define the left\--continuous derivator $g:[0,1]\to\mathbb{R}$ as
 \begin{equation}\label{eq:example:g}
 g(t):=\begin{dcases}
 \int_{x_0}^t h(s) \operatorname{d} s , & \mbox{if } t\in[x_0,1],\\
 0, & \mbox{if } t\in [0,x_0].
 \end{dcases}
 \end{equation}
 The derivator $g$ is absolutely continuous on $[x_0,1]$, thus, it follows from~\cite[Proposition~5.3]{PR}, that $g$ is of bounded variation. Furthermore, $\widetilde g(x)=x-x_0$ for $x\in[x_0,1]$. We will now study the values of the function
\begin{equation*}
\psi(t):=\frac{g(t)-g(x_0)}{\widetilde{g}(t)-\widetilde{g}(x_0)}=\frac{g(t)}{\widetilde{g}(t)},\quad t\in(x_0,1],
\end{equation*}
on the points of the sequence $(x_n)_{n\in{\mathbb N}}$. Using the definition of $g$ and the fact that $h$ is constant on each interval of the form $(x_{k+1},x_k]$,
\begin{equation}\label{eq2n}
\psi(x_{2n})=\frac{\sum_{k=n}^\infty(x_{2k}-x_{2k+1})-\sum_{k=n}^\infty(x_{2k+1}-x_{2k+2})}{\sum_{k=n}^\infty(x_{2k}-x_{2k+1})+\sum_{k=n}^\infty(x_{2k+1}-x_{2k+2})}=\frac{x_{2n}+2\sum_{k=n}^\infty (x_{2k+2}-x_{2k+1})}{x_{2n}}=\frac{x_{2n}+2s_n}{x_{2n}},
\end{equation}
where $s_n:=\sum_{k=n}^\infty (x_{2k+2}-x_{2k+1})$. On the other hand,
\begin{equation}
\label{eq2n1}\psi(x_{2n+1})=\frac{x_{2n}+2s_n-(x_{2n}-x_{2n+1})}{x_{2n+1}}=\frac{x_{2n+1}+2s_n}{x_{2n+1}}.
\end{equation}
Thus, solving for $2s_n$ in~\eqref{eq2n} and~\eqref{eq2n1} and equating,
\[
(\psi(x_{2n})-1)x_{2n}=(\psi(x_{2n+1})-1)x_{2n+1}=(\psi(x_{2n+1})-1)(1-\alpha_n)x_{2n},
\]
we conclude that
\begin{equation*}
\psi(x_{2n+1})=\frac{\psi(x_{2n})-1}{1-\alpha_n}+1.
\end{equation*}
If we can prove that for every $n\in\mathbb{N}$, $\psi(x_{2n})=\alpha_n$, then we would have that $\psi(x_{2n+1})=0$.

On the other hand, evaluating~\eqref{eq2n1} on $n-1$,
\begin{equation}\label{eq2n2}
\psi(x_{2n-1})=\frac{x_{2n-1}+2s_{n-1}}{x_{2n-1}}=\frac{x_{2n-1}+2(s_{n}+x_{2n}-x_{2n-1})}{x_{2n-1}}=\frac{-x_{2n-1}+2x_{2n}+2s_{n}}{x_{2n-1}}.
\end{equation}
Solving for $2s_n$ in~\eqref{eq2n} and~\eqref{eq2n2} and equating,
\begin{equation*}
(\psi(x_{2n-1})+1)x_{2n-1}-2x_{2n}=(\psi(x_{2n})-1)x_{2n},
\end{equation*}
so
\begin{equation*}
(\psi(x_{2n-1})+1)x_{2n-1}=(\psi(x_{2n})+1)x_{2n}=\frac{\psi(x_{2n})+1}{1+\alpha_n}x_{2n-1}.
\end{equation*}
Thus,
\begin{equation*}
\psi(x_{2n})=(\psi(x_{2n-1})+1)(1+\alpha_{n})-1.
\end{equation*}
If we can prove that for every $n\in\mathbb{N}$, $\psi(x_{2n-1})=0$, then we would have that $\psi(x_{2n})=\alpha_n$.

Using~\eqref{eq2n2} and~\eqref{eq:sum for psi(x_1)}, we obtain
\begin{align*}
\psi(x_{1})= & \frac{-x_{1}+2x_2+2s_1}{x_{1}}=\frac{-x_{1}+2x_2+2\sum_{k=1}^\infty (x_{2k+2}-x_{2k+1})}{x_{1}}\\
= & \frac{-x_{1}+2x_2+2\sum_{k=1}^\infty\left( \frac{x_1}{1+\alpha_{k+1}}\prod_{j=1}^k\frac{1-\alpha_j}{1+\alpha_j}-x_{1}\prod_{j=1}^k\frac{1-\alpha_j}{1+\alpha_j}\right) }{x_{1}}\\
= &-1+\frac{2}{1+\alpha_1}+2\sum_{k=1}^\infty\left( \frac{1}{1+\alpha_{k+1}}-1\right) \prod_{j=1}^k\frac{1-\alpha_j}{1+\alpha_j}\\
= & \frac{1-\alpha_1}{1+\alpha_1}-2\sum_{k=1}^\infty\frac{\alpha_{k+1}}{1+\alpha_{k+1}}\prod_{j=1}^k\frac{1-\alpha_j}{1+\alpha_j}\\
= & 0.
\end{align*}
Hence, it follows that $\psi(x_{2n-1})=0$ for all $n\in\mathbb{N}$.

Since $\psi(x_{2n-1}) =0$ and $x_{2n-1}\to x_0^+$, then \[ \varphi(x_0)=\liminf_{t\to x_0^+}\left|\psi(t)\right|\leqslant \liminf_{n\to\infty}\psi(x_{2n-1}) =0.
\]
\end{exa}

Now consider a function $f\in{\mathcal C}([x_0,1],{\mathbb R})$ satisfying the following properties for every ${n\in{\mathbb N}}$:
\begin{enumerate}
	\item $f(x_n)=f(x_0)=0$,
	\item $f(x)< 0$ if $x\in (x_{2n},x_{2n-1})$ and
	$f(x)> 0$ if $x\in (x_{2n+1},x_{2n})$,
 \item $|f(x)|\leqslant \frac{x_{n}^{1+r}-x_{n+1}^{1+r}}{x_n-x_{n+1}}$ if $x\in[x_{n+1},x_n]$, for some $r\in{\mathbb R}^+$ such that $\lim\limits_{n\to \infty }x_n^{r}= 0$, and $\lim\limits_{n\to \infty }\alpha_nx_{2n}^{-r}=0$ if $x_0=0$.
	\item $\int_{x_{n+1}}^{x_n}|f(s)|\operatorname{d} s=\frac{1}{2}(x_{n}^{1+r}-x_{n+1}^{1+r})$.
\end{enumerate}

Observe that, applying the mean value theorem to the function $\eta(x):=x^{1+r}$ on $[x_{n+1},x_n]$, we conclude using condition~3 that for all $x\in[x_{n+1},x_n]$, we have $|f(x)|\leqslant (1+r)c_n^{r}$ for some $c_n\in(x_{n+1},x_n)$.
Since $0\leqslant c_n^{r}\leqslant x_n^{r}$ and $x_n^{r}\xrightarrow{n\to\infty}0$, we have that $\lim_{x\to x_0^+}|f(x)|=0$, which guarantees continuity at~$x_0$.

A simple example of such a function is one that, on each interval $[x_{n+1},x_n]$ its graph has the shape of a triangle with two vertices at the points $(x_{n+1},0)$ and $(x_{n},0)$ and its third vertex at the point
\[  \left( \frac{1}{2}(x_n+x_{n+1}), (-1)^{n}\frac{x_{n}^{1+r}-x_{n+1}^{1+r}}{x_n-x_{n+1}}\right).\]
Given that $f$ is continuous and bounded we have that $f$ is $\widetilde{g}$\--continuous on $[x_0,1]$ and, therefore, $f \in {\mathcal B}{\mathcal C}_{g}([x_0,1],\mathbb{R})$, so let, for $t\in[x_0,1]$,
\[  F(t):=\int_{[x_0,t)}f(s)\operatorname{d}\mu_g(s)=\int_{x_0}^t|f(s)|\operatorname{d} s,\]
where the equality holds because $|g'(t)|=1$ for a.e. $t\in(x_0,1]$ and $f$ is chosen such that $f(t)g'(t)=|f(t)|$ for a.e. $t\in(x_0,1]$. Observe that $f$ is non\--negative where $g$ is nondecreasing and non\--positive where $g$ is decreasing.
In particular,
\[  F(x_n)=\sum_{k=n}^\infty\frac{1}{2}(x_{k}^{1+r}-x_{k+1}^{1+r})=\frac{1}{2}(x_n^{1+r}-x_0).
\]
If $F'_g(x_0)$ exists, it has to be computed as
\begin{align*}
F'_g(x_0)=& \lim_{t\to x_0^+}\frac{F(t)-F(x_0)}{g(t)-g(x_0)}=\lim_{n\to \infty }\frac{F(x_n)}{g(x_n)}=\lim_{n\to \infty }\frac{\frac{1}{2}(x_n^{1+r}-x_0)}{g(x_n)}=\frac{1}{2}\lim_{n\to \infty }\frac{\widetilde g(x_n)}{g(x_n)}\frac{x_n^{1+r}-x_0}{x_n-x_0}\\ = &\frac{1}{2}\lim_{n\to \infty }\frac{1}{\psi(x_n)}\frac{x_n^{1+r}-x_0}{x_n-x_0}.
\end{align*}
We have two cases. If $x_0=0$, then
\[  F'_g(x_0)=\frac{1}{2}\lim_{n\to \infty }\frac{x_n^{r}}{\psi(x_n)}=\frac{1}{2}\lim_{n\to \infty }\frac{x_{2n}^{r}}{\psi(x_{2n})}=\frac{1}{2}\lim_{n\to \infty }\frac{x_{2n}^{r}}{\alpha_n}=\infty,\]
and we arrive to a contradiction.

If $x_0\ne 0$, then
\[
F'_g(x_0)= -\frac{x_0}{2}\lim_{n\to \infty }\frac{1}{(x_n-x_0)\psi(x_n)} =-\infty,
\]
given that $(x_n-x_0)\psi(x_n)\xrightarrow{n\to\infty} 0$, as the non\--negative sequence $\left( \psi(x_n)\right) _{n\in\mathbb{N}}$ is bounded and $x_n-x_0 \xrightarrow{n\to\infty}0^+$. Thus, we arrive to a contradiction as well.

\begin{exa}\label{exa2}In this example, we choose a specific set of parameters in Example~\ref{exa1} in order to show that the restrictions imposed on those parameters can be met.
Let us consider the sequence $(\alpha_n)_{n\in\mathbb N}$ defined as
\[
\begin{dcases}
\alpha_1=\frac{1}{2}, & \\
\alpha_n=\frac{1}{n}, & \mbox{for all $n\in{\mathbb N}$, $n\geqslant2$}.
\end{dcases}
\]
Observe that~\eqref{eq:sum for psi(x_1)} holds. Indeed,
\begin{align*}
\sum_{k=1}^\infty\frac{\alpha_{k+1}}{1+\alpha_{k+1}}\prod_{j=1}^k\frac{1-\alpha_j}{1+\alpha_j}
=&\frac{\alpha_2}{1+\alpha_2}\frac{ 1-\alpha_1}{1+\alpha_1}+\frac{ 1-\alpha_1}{1+\alpha_1} \sum_{k=2}^\infty\frac{\alpha_{k+1}}{1+\alpha_{k+1}}\prod_{j=2}^k\frac{1-\alpha_j}{1+\alpha_j}\\
=&\frac{ 1-\alpha_1}{1+\alpha_1}\left(\frac{\alpha_2}{1+\alpha_2}+ \sum_{k=2}^\infty\frac{\frac{1}{k+1}}{1+\frac{1}{k+1}}\prod_{j=2}^k\frac{1-1/j}{1+1/j}\right)\\
=&\frac{ 1-1/2}{1+1/2}\left(\frac{1/2}{1+1/2}+ \sum_{k=2}^\infty\frac{1}{k+2}\prod_{j=2}^k\frac{j-1}{j+1}\right)\\
=&\frac{1}{3}\left(\frac{1}{3}+ \sum_{k=2}^\infty\frac{2}{k+2}\frac{(k-1)!}{(k+1)!}\right)\\
=&\frac{1}{3}\left(\frac{1}{3}+2 \sum_{k=2}^\infty\frac{1}{(k+2)(k+1)k}\right)\\
=&{\frac{1}{3}\left(\frac{1}{3}+ \sum_{k=2}^\infty \left(\frac{1}{k(k+1)} -\frac{1}{(k+1)(k+2)}\right)\right)}\\
=&{\frac{1}{3}\left(\frac{1}{3}+\frac{1}{2\cdot3}\right)=\frac{1}{6}}\\
=& \frac{1}{2}\frac{ 1-1/2}{1+1/2}=\frac{1}{2}\frac{ 1-\alpha_1}{1+\alpha_1}.
\end{align*}

We now define the sequence $(x_n)_{n\in{\mathbb N}}\subset [0,1]$ as follows:
\begin{equation}\label{eq:x_n}
x_1=1,\quad x_{2n}= \frac{1}{1+\alpha_n}x_{2n-1},\quad x_{2n+1}=(1-\alpha_n)x_{2n},\quad {n\in{\mathbb N}}.
\end{equation}
Substituting the values of~$(\alpha_n)_{n\in{\mathbb N}}$, we obtain
\[
x_1=1,\quad x_2=\frac{2}{3},\quad x_3=\frac{1}{3}, \quad x_4=\frac{2}{9},\quad x_5=\frac{1}{9}, \quad x_6=\frac{1}{12},\quad x_7=\frac{1}{18},\quad x_8=\frac{2}{45}\cdots
\]
In particular, for every $n\geqslant 2$
\[
 x_{2n}=\frac{1}{1+\alpha_n}x_{2n-1}=\frac{n}{n+1}x_{2n-1},\quad x_{2n+1}=(1-\alpha_n)x_{2n}=\frac{n-1}{n}x_{2n},\quad {n\in{\mathbb N}}.
\]
Therefore, by induction, for every $n\geqslant 2$, we obtain
\[
x_{2n}=\frac{2}{3(n-1)(n+1)},
\quad\text{ and }\quad
x_{2n+1}=\frac{2}{3n(n+1)}.
\]
Observe that $x_n\xrightarrow{n\to\infty}x_0=0$. In Figure~\ref{fig:graph of g}, we illustrate the graph of the derivator $g:[0,1]\to\mathbb{R}$ defined in~\eqref{eq:example:g} associated to the sequence $(x_n)_{n\in\mathbb{N}}$.

\begin{figure}[h]
 \centering
 \includegraphics[width=17cm]{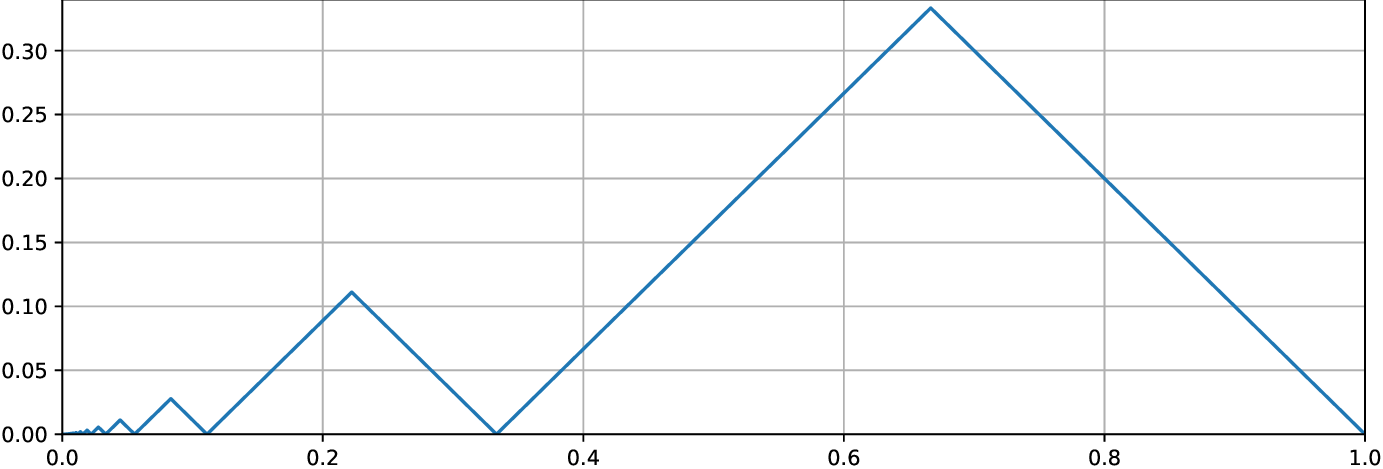}
 \caption{The graph of the derivator $g$ associated to the sequence~\eqref{eq:x_n}.}\label{fig:graph of g}
\end{figure}

Now, let $
r=\frac{1}{3}$.
Observe that
\[
\lim_{n\to\infty}x_{2n}^{r}= \lim_{n\to\infty} \left(\frac{2}{3(n-1)(n+1)}\right)^{\frac{1}{3}}= \lim_{n\to\infty} \left(\frac{2 }{3(n^2-1)}\right)^{\frac{1}{3}}=0,
\]
\[
\lim_{n\to\infty}x_{2n+1}^{r}= \lim_{n\to\infty} \left(\frac{2}{3n(n+1)}\right)^{\frac{1}{3}}= \lim_{n\to\infty} \left(\frac{2 }{3(n^2+n)}\right)^{\frac{1}{3}}=0,
\]
and
\[
\lim_{n\to\infty}\frac{\alpha_n}{x_{2n}^{r}}=\lim_{n\to\infty}\frac{1}{n\left(\frac{2}{3(n-1)(n+1)}\right)^{\frac{1}{3}}}=
\lim_{n\to\infty}\frac{1}{\left(\frac{2n^3 }{3(n^2-1)}\right)^{\frac{1}{3}}}=0.
\]
Figure~\ref{fig:graph of f} shows the graph of the function $f: [0,1]\to{\mathbb R}$ whose graph on each interval $[x_{n+1},x_n]$ takes the shape of a triangle, with two vertices at the points $(x_{n+1},0)$ and $(x_{n},0)$ and its third vertex at the point
\[  \left( \frac{1}{2}(x_n+x_{n+1}), (-1)^{n}\frac{x_{n}^{1+r}-x_{n+1}^{1+r}}{x_n-x_{n+1}}\right).\]

\begin{figure}[ht]
\includegraphics[width=.49\textwidth]{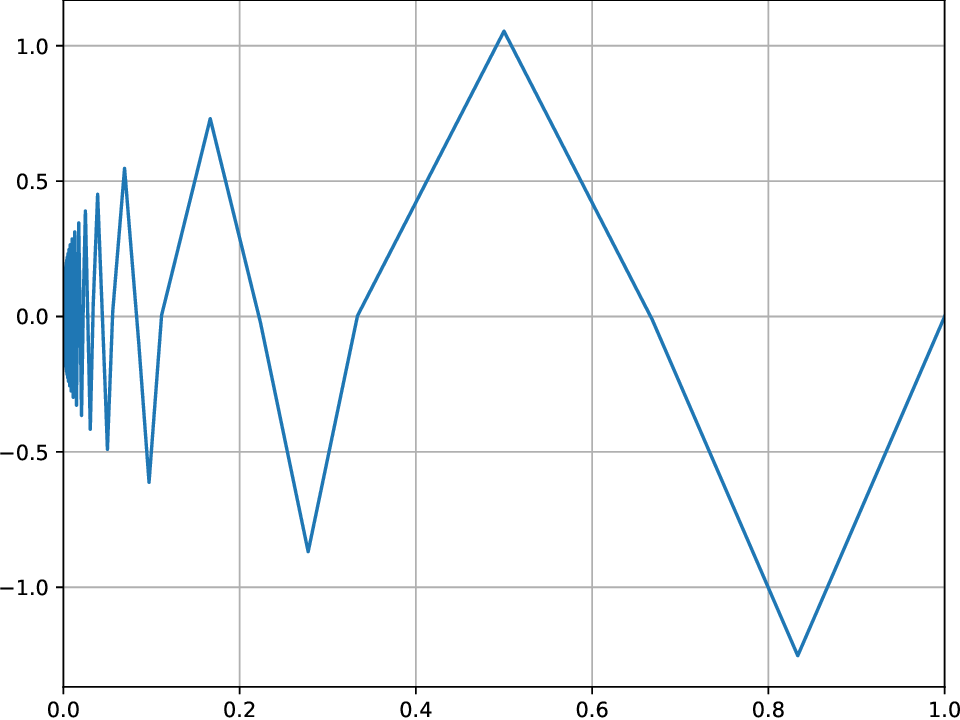}
 \includegraphics[width=.49\textwidth]{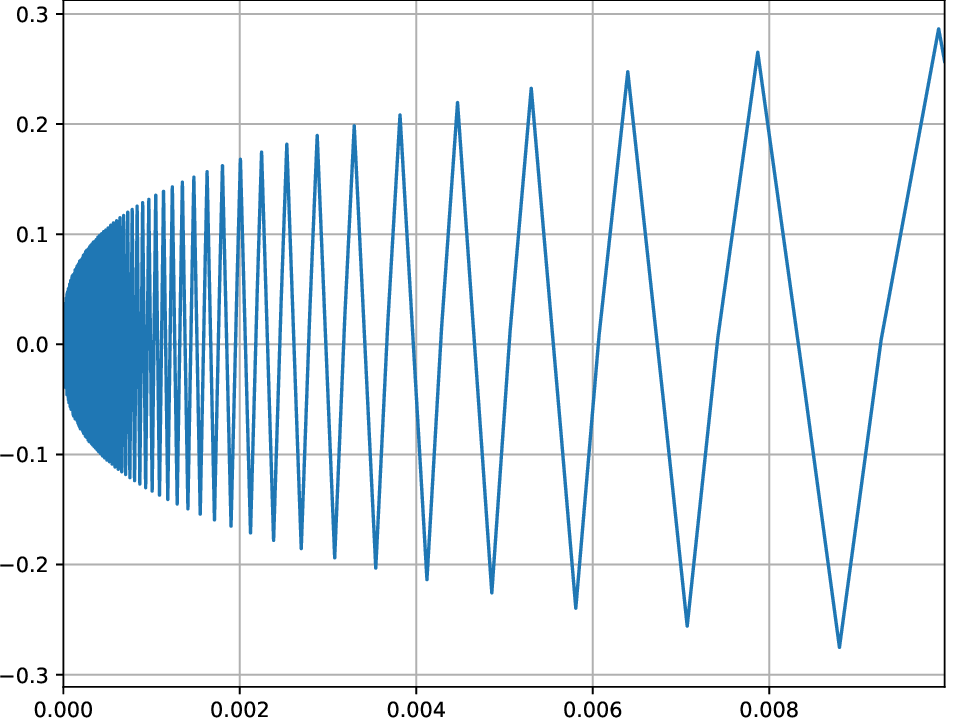}
 \caption{Graph of f (left) and zoomed-in view (right).}\label{fig:graph of f}
\end{figure}

Now, let us compute $F(t)$ explicitly: Let $t\in[x_{k+1},x_k]$ for some $k\in\mathbb{N}$.
\begin{align*}
 F(t)= & \int_0^t |f|\,\operatorname{d}s=\int_0^1 |f|\,\operatorname{d}s-\int_t^1 |f|\,\operatorname{d}s=F(1)-\int_t^1 |f|\,\operatorname{d}s \\
 = & \frac{1}{2}(x_1^{1+r}-x_0)- \int_t^1 |f|\,\operatorname{d}s= \frac{1}{2} -\left(\sum_{j=1}^{k-1}\int_{x_{j+1}}^{x_j} |f|\,\operatorname{d}s + \int_{t}^{x_k} |f|\,\operatorname{d}s\right).
\end{align*}
Let us set $x_{\text{mid},j}:=\frac{x_j+x_{j+1}}{2}$. Then,
\begin{align*}
 F(t) = & \frac{1}{2} -\left(\sum_{j=1}^{k-1}2\int_{x_{j+1}}^{x_{\text{mid},j}} |f|\,\operatorname{d}s + \int_{t}^{x_k} |f|\,\operatorname{d}s\right)\\
 = & \frac{1}{2} -\left(\sum_{j=1}^{k-1}2\int_{x_{j+1}}^{x_{\text{mid},j}} \frac{x_{j}^{1+r}-x_{j+1}^{1+r}}{x_j-x_{j+1}}\frac{s-x_{j+1}}{x_{\text{mid},j}-x_{j+1}}\,\operatorname{d}s + \int_{t}^{x_k} |f|\,\operatorname{d}s\right)\\
 = & \frac{1}{2} -\left(\sum_{j=1}^{k-1}\frac{x_{j}^{1+r}-x_{j+1}^{1+r}}{x_j-x_{j+1}}(x_{\text{mid},j}-x_{j+1}) + \int_{t}^{x_k} |f|\,\operatorname{d}s\right)\\
 = & \frac{1}{2} -\left(\sum_{j=1}^{k-1}\frac{x_{j}^{1+r}-x_{j+1}^{1+r}}{x_j-x_{j+1}}\frac{x_{j}-x_{j+1}}{2} + \int_{t}^{x_k} |f|\,\operatorname{d}s\right)\\
 = & \frac{1}{2} -\left(\sum_{j=1}^{k-1}\frac{x_{j}^{1+r}-x_{j+1}^{1+r}}{2} + \int_{t}^{x_k} |f|\,\operatorname{d}s\right)\\
 = & \frac{1}{2} -\left(\frac{1}{2}\sum_{j=1}^{k-1}(x_{j}^{4/3}-x_{j+1}^{4/3}) + \int_{t}^{x_k} |f|\,\operatorname{d}s\right)\\
 = & \frac{1}{2} -\left(\frac{1}{2}(x_{1}^{4/3}-x_{k}^{4/3}) + \int_{t}^{x_k} |f|\,\operatorname{d}s\right)= \frac{1}{2}x_{k}^{4/3}- \int_{t}^{x_k} |f|\,\operatorname{d}s.\\
\end{align*}
Now, let us compute $\int_{t}^{x_k} |f|\,\operatorname{d}s$. We distinguish two cases:

\textbf{Case 1: } If $t\in[x_{k+1},x_{\text{mid},k}]$, then
\begin{align*}
 \int_{t}^{x_k} |f|\,\operatorname{d}s= & \int_{t}^{x_{\text{mid},k}} |f|\,\operatorname{d}s + \int_{x_{\text{mid},k}}^{x_k} |f|\,\operatorname{d}s\\
 = & \int_{t}^{x_{\text{mid},k}}\frac{x_{k}^{1+r}-x_{k+1}^{1+r}}{x_k-x_{k+1}}\frac{s-x_{k+1}}{x_{\text{mid},k}-x_{k+1}} \,\operatorname{d}s + \int_{x_{\text{mid},k}}^{x_k} \frac{x_{k}^{1+r}-x_{k+1}^{1+r}}{x_k-x_{k+1}}\frac{x_k-s}{x_{k}-x_{\text{mid},k}} \,\operatorname{d}s \\
 = & \frac{x_{k}^{1+r}-x_{k+1}^{1+r}}{x_k-x_{k+1}} \left(\int_{t}^{x_{\text{mid},k}}\frac{s-x_{k+1}}{x_{\text{mid},k}-x_{k+1}} \,\operatorname{d}s + \int_{x_{\text{mid},k}}^{x_k} \frac{x_k-s}{x_{k}-x_{\text{mid},k}} \,\operatorname{d}s\right) \\
 = & \frac{1}{2}\frac{x_{k}^{1+r}-x_{k+1}^{1+r}}{x_k-x_{k+1}} \left(\frac{(x_{\text{mid},k}-x_{k+1})^2-(t-x_{k+1})^2}{x_{\text{mid},k}-x_{k+1}} + x_k-x_{\text{mid},k}\right) \\
 = & \frac{x_{k}^{1+r}-x_{k+1}^{1+r}}{(x_k-x_{k+1})^2} \left((x_{\text{mid},k}-x_{k+1})^2-(t-x_{k+1})^2 + (x_{\text{mid},k}-x_{k+1})^2\right) \\
 =& \frac{x_{k}^{1+r}-x_{k+1}^{1+r}}{(x_k-x_{k+1})^2} \left(\frac{(x_k-x_{k+1})^2}{2}-(t-x_{k+1})^2\right)\\
 =& \frac{x_{k}^{4/3}-x_{k+1}^{4/3}}{(x_k-x_{k+1})^2} \left(\frac{(x_k-x_{k+1})^2}{2}-(t-x_{k+1})^2\right).
\end{align*}

\textbf{Case 2: } If $t\in[x_{\text{mid},k},x_{k}]$, then
\[
 \int_{t}^{x_k} |f|\,\operatorname{d}s= \int_{t}^{x_k} \frac{x_{k}^{1+r}-x_{k+1}^{1+r}}{x_k-x_{k+1}}\frac{x_k-s}{x_{k}-x_{\text{mid},k}} \,\operatorname{d}s = \frac{x_{k}^{4/3}-x_{k+1}^{4/3}}{(x_k-x_{k+1})^2}(x_k-t)^2 .
\]

In Figure~\ref{fig:graph of g and F}, we illustrate the behavior of the function~$F$ and the derivator~$g$ around the point $x_0=0$.
Additionally, in Figure~\ref{fig:graph of Q)}, we depict the behavior of the function $Q:=\frac{F(\cdot)-F(x_0)}{g(\cdot)-g(x_0)}=\frac{F}{g}$ which is defined on $(0,1]\setminus \{x_{2n-1}\}_{n\in\mathbb{N}}$. The behavior of~$Q$ illustrates the lack of $g$\--differentiability of~$F$ at~$x_0=0$.
\begin{figure}[ht]
\begin{subfigure}{0.49\textwidth}
\centering
 \includegraphics[width=.99\linewidth]{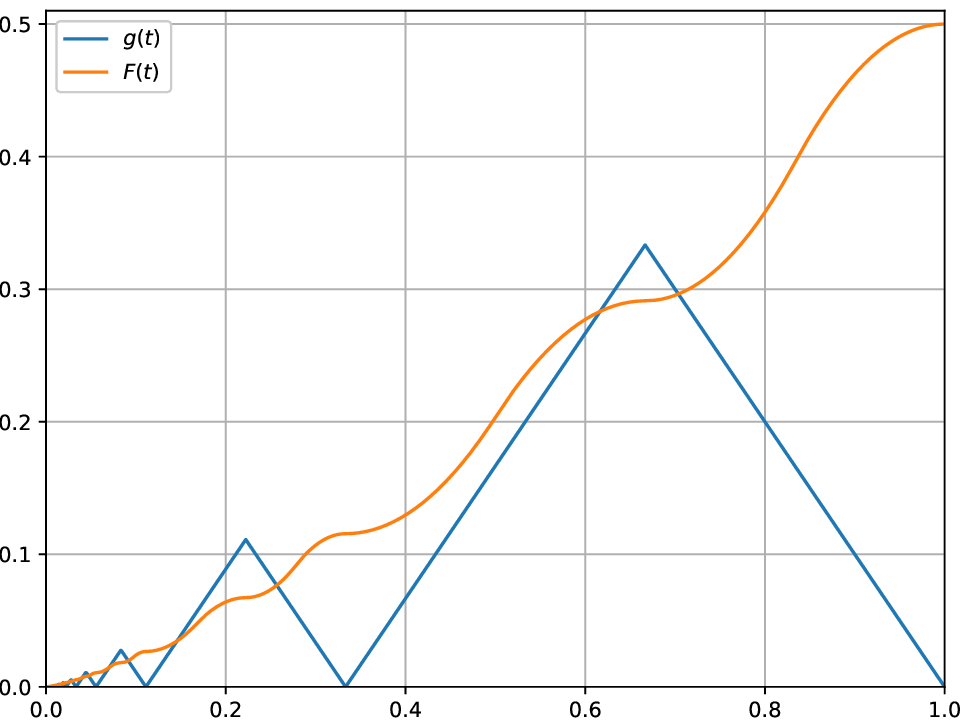}
\caption{Graph of $F$ and $g$.}\label{fig:graph of g and F}
\end{subfigure}
\begin{subfigure}{0.49\textwidth}
\centering
\includegraphics[width=.99\linewidth]{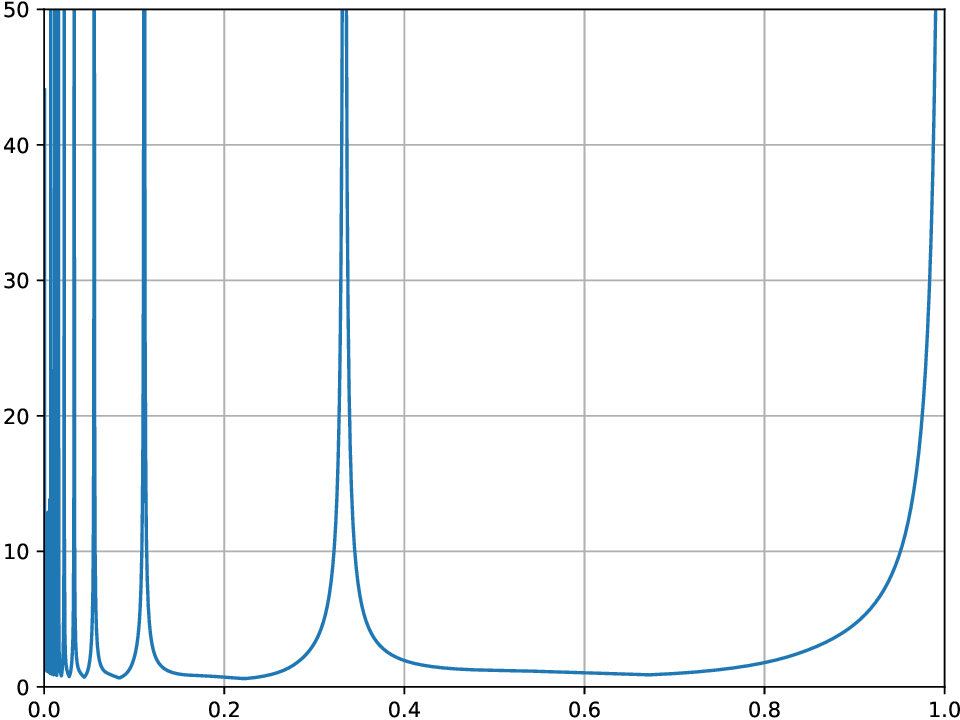}
\caption{Graph of $Q$.}\label{fig:graph of Q)}
\end{subfigure}\caption{Illustration of the fact that $F$ cannot have a $g$\--derivative at $0$.}
\end{figure}
\end{exa}
\begin{rem}
 In Example~\ref{exa2}, the function $F$ is not $g$\--differentiable at $x_0=0$. Nevertheless, since $\widetilde{g}=\operatorname{id}_{[0,1]}$, and $f$ is continuous on $[0,1]$ (and in particular Lebesgue integrable), it follows that $F$ is differentiable at $x_0=0$ in the classical sense. This implies, in particular, that $\widetilde{g}$\--differentiability does not necessarily yield $g$\--differentiability.
\end{rem}

We will now show that the problems occurring in Examples~\ref{exa1} and~\ref{exa2} happen for any derivator $g$, which means that Theorem~\ref{thmfcf} is optimal in its assumptions. In order to illustrate this, we will need some auxiliary results.

\FloatBarrier
\begin{thm}[{\cite[Theorem 2.3]{MT}}]\label{teoroutermeasureg}
	Let $g:\mathbb R\to\mathbb R$ be a nondecreasing and left\--continuous function and $\mu_g^*$ be the exterior measure associated to $g$. Then, for any $A\in\mathcal P(\mathbb R)$,
	\begin{equation*}
		\mu^*_g(A)=\inf\left\{\sum_{n=1}^\infty (g(b_n)-g(a_n)): A\subset \bigcup_{n=1}^\infty [a_n,b_n),\ \{[a_n,b_n)\}_{n=1}^\infty\subset \mathcal C\text{ pairwise disjoint}\right\},
	\end{equation*}
with $\mathcal C=\{[a,b):a,b\in\mathbb{R},\, a<b\}$.
\end{thm}
\FloatBarrier
The following result is a density theorem that will serve as a key tool for approximating $g$\--integrable functions by $g$\--continuous functions. For the sake of generality, we present the result in the most general setting, where $J\subset\mathbb{R}$ is the definition domain of $g$ instead of restricting on $[a,b]$, and we consider $g:J\to{\mathbb R}$ with no unbounded constancy interval of the form $(-\infty,x]$. In doing so, we retain the notations introduced in previous sections but with the understanding that $g$ is defined on $J$ rather than on the interval $[a,b]$.

For simplicity, given $n\in{\mathbb N}$ and $A=\{(x_k,y_k)\}_{k=1}^n\subset{\mathbb R}^2$, such that $x_j<x_k$ if $j<k$, we will denote by $P_A(x)$ the piecewise linear function
\[ P_A(x):=\begin{dcases}y_1, & \text{ if $x\in(-\infty,{x_1}]$},\\
	y_k+\frac{y_{k+1}-y_k}{x_{k+1}-x_k}(x-x_k), & \text{ if } x\in({x_k},{x_{k+1}}),\ k\in\{1,\dots,n-1\},\\
	y_n, &\text{ if $x\in ({x_n},\infty)$}.
	\end{dcases}\]
Observe that $\{(x_1,y_1),(x_1,y_1),(x_2,y_2)\}=\{(x_1,y_1),(x_2,y_2)\}$, so there is no problem if, when defining the set $A$, there are repeated points, as long as there are no two points with the same $x$ coordinate and different $y$ coordinates.

\begin{thm}\label{density}
Let $I\subset J\subset{\mathbb R}$ be two fixed intervals, and $g:J\to{\mathbb R}$ a nondecreasing and left-continuous derivator.
	\begin{enumerate}
\item	${\mathcal B}{\mathcal C}_g(I,{\mathbb R})\cap L^1_g(I,{\mathbb R})$ is dense in $L^1_g(I,{\mathbb R})$.
\item If $f:I\to[c,d]$ is a $g$\--integrable function then, for every $\varepsilon\in{\mathbb R}^+$, there exists $h\in{\mathcal B}{\mathcal C}_g({I},[c,d])$ such that $\int_{I}|f-h|\operatorname{d}\mu_g<\varepsilon$.
\item If $J=I=[a,b]$ for some $a,b\in\mathbb{R}$ such that $a^*\not\in D_g$, $g(a)<g(b)$, $\alpha,\beta\in[c,d]$, and $f:I\to[c,d]$ is a $g$\--integrable function, then, for every $\varepsilon>0$, there exists $h\in{\mathcal B}{\mathcal C}_g(I,[c,d])$ such that, $h(a)=\alpha$, $h(b)=\beta$ and $\int_{[a,b)}|f-h|\operatorname{d}\mu_g<\varepsilon$.
\item If $J=I=[a,b]$ for some $a,b\in\mathbb{R}$ such that $a^*\in D_g$, $\beta\in[c,d]$, and $f:I\to[c,d]$ is a $g$\--integrable function, then, for every $\varepsilon>0$, there exists $h\in{\mathcal B}{\mathcal C}_g(I,[c,d])$ such that $h(a)=f(a^*)$, $h(b)=\beta$ and $\int_{[a,b)}|f-h|\operatorname{d}\mu_g<\varepsilon$.
\end{enumerate}
	\end{thm}
The proof of points~1, and~2 of this theorem will be given in the case where $J=\mathbb{R}$, the proof for other cases of $J$ follows analogously.
\begin{proof}
1. We prove the case where $I={\mathbb R}$ (the rest would be analogous). Let \[  C := \overline{{\mathcal B}{\mathcal C}_g(I,{\mathbb R})\cap L^1_g(I,{\mathbb R})}.\]  $C$ is a closed vector subspace of $L^1_g(I,{\mathbb R})$. We now show, step by step, that $C=L^1_g(I,{\mathbb R})$.

We will be considering the function $g^\dagger(y):=\inf \{t \in {\mathbb R}: g(t) \geqslant y\},\ y \in g({\mathbb R})$\---for more information on the properties of this function see~\cite{Tojo2025connectingSDEnODE}.

	 \emph{$\diamond$ Step 1: $\chi_{[a,b)}\in C$ with $a, b \in \mathbb{R}$, $a < b$.} We start by considering the case $a<g^{\dagger}(g(b))$. Fix $\varepsilon>0$. Since $g$ is left\--continuous with no unbounded constancy intervals, we can take $t=g^{\dagger}(g(a))$ if $g(g^{\dagger}(g(a)))<g(a)$ and, otherwise, $t<a$ such that
\[
g(t)<g(a) \text{ with } g(a)-g(t)<\frac{\varepsilon}{2};
\]
and $s=g^\dagger(g(b))$ if $g(g^\dagger(g(b)))<g(b)$ and, otherwise, $s\in [a,b)$ such that
and
\[
g(s)<g(b) \text{ with } g(b)-g(s)<\frac{\varepsilon}{2}.
\]

Let $f=P_A$ where \[ A=\{(g(t),0),(g(a),1),(g(s),1),(g(b),0)\},\]
an example is shown in Figure~\ref{fig:proof:f:step1:ex1}.

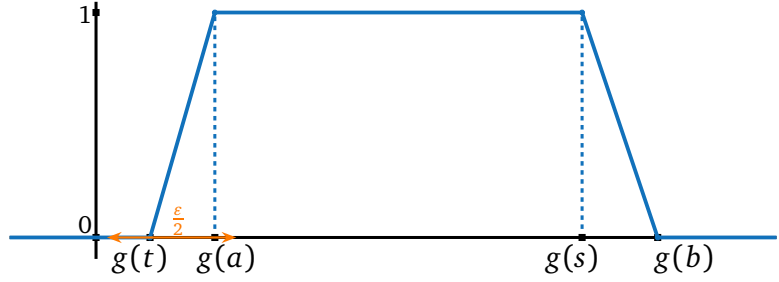
\begin{figure}[h]
	\centering
	\resizebox{0.6\textwidth}{!}{%
		\begin{circuitikz}
			\tikzstyle{every node}=[font=\LARGE]
			\draw [ line width=2pt](-2.75,3) to[short] (15,3);
			\draw [line width=2pt, short] (-0.75,8.5) -- (-0.75,2.5);

			\draw [ color={rgb,255:red,17; green,113; blue,187}, line width=2pt, dashed] (2,3.2) -- (2,8.25);
			\draw [ color={rgb,255:red,17; green,113; blue,187}, line width=2pt, dashed] (10.5,3.2) -- (10.5,8.25);

			\node at (-0.75,3) [squarepole,line width=1.5pt] {};
			\node [font=\Large] at (-1,3.3) {$0$};
			\node at (-0.75,8.25) [squarepole,line width=1.5pt] {};
			\node [font=\Large] at (-1,8.2) {$1$};

			\node at (0.5,3) [squarepole,line width=1.5pt] {};
			\node at (0.5,3) [circ, color={rgb,255:red,17; green,113; blue,187},line width=0.001pt] {};
			\node [font=\LARGE] at (0.25,2.5) {$g(t)$};

			\node at (2,3) [squarepole,line width=1.5pt] {};
			\node at (2,8.25) [circ, color={rgb,255:red,17; green,113; blue,187},line width=0.001pt] {};
			\node [font=\LARGE] at (2.25,2.5) {$g(a)$};


			\node at (10.5,3) [squarepole,line width=1.5pt] {};
			\node at (10.5,8.25) [circ, color={rgb,255:red,17; green,113; blue,187},line width=0.001pt] {};
			\node [font=\LARGE] at (10.25,2.5) {$g(s)$};

			\node at (12.25,3) [squarepole,line width=1.5pt] {};
			\node at (12.25,3) [circ, color={rgb,255:red,17; green,113; blue,187},line width=0.001pt] {};
			\node [font=\LARGE] at (12.85,2.5) {$g(b)$};

			\draw [ color={rgb,255:red,17; green,113; blue,187}, line width=2.5pt] (-2.75,3) -- (0.5,3);
			\draw [ color={rgb,255:red,17; green,113; blue,187}, line width=2.5pt] (2,8.25) -- (0.5,3);
			\draw [ color={rgb,255:red,17; green,113; blue,187}, line width=2.5pt] (2,8.25) -- (10.5,8.25);
			\draw [ color={rgb,255:red,17; green,113; blue,187}, line width=2.5pt] (10.5,8.25) -- (12.25,3);
			\draw [ color={rgb,255:red,17; green,113; blue,187}, line width=2.5pt] (12.25,3) -- (15,3);

			\draw [ color={rgb,255:red,254; green,121; blue,0}, line width=1.5pt, <->, >=Stealth] (-0.5,3) -- (2.5,3);
			\node [font=\Large, color={rgb,255:red,254; green,121; blue,0}] at (1.15,3.4) {$\frac{\varepsilon}{2}$};

		\end{circuitikz}
	}%
	\caption{The continuous function $f$ with $g^{\dagger}(g(a))<a$, $g^{\dagger}(g(b))=b$, $g(a)-g(t)<\frac{\varepsilon}{2}$.}\label{fig:proof:f:step1:ex1}
\end{figure}

 Since $f$ is continuous, $f\circ g$ is $g$\--continuous and of compact support (with respect to the usual topology of ${\mathbb R}$), so $f\circ g\in L^1_g(I,{\mathbb R})$. Observe that $f\circ g=\chi_{[a,b)}$ on the set $(-\infty,t]\cup[a,s]\cup[b,\infty)$. Moreover, for $u\in\{a,b\}$, if $g^\dagger(g(u))<u$ then $f\circ g(u)=\chi_{[a,b)}(u)$, and $(g^\dagger(g(u)),u)\subset C_g$. Therefore,
\FloatBarrier
\begin{align*}\|f\circ g-\chi_{[a,b)}\|_{L_g^1(\mathbb{R})}
= & \int_{(t,g^{\dagger}(g(a)))}|f\circ g-\chi_{[a,b)}|\operatorname{d} \mu_g+ \int_{(s,g^{\dagger}(g(b)))}|f\circ g-\chi_{[a,b)}|\operatorname{d}\mu_g\\
\leqslant & \int_{(t,g^{\dagger}(g(a)))}(1-f\circ g) \operatorname{d} \mu_g+\int_{(s,g^{\dagger}(g(b)))}(1-f\circ g)\operatorname{d}\mu_g\\
\leqslant & \int_{[t,g^{\dagger}(g(a)))}1 \operatorname{d} \mu_g+\int_{[s,g^{\dagger}(g(b)))} 1\operatorname{d}\mu_g\\
 = & g\big(g^{\dagger}(g(a))\big)-g(t)+g\big(g^{\dagger}(g(b))\big)-g(s).
\end{align*}
\FloatBarrier
Given that either $t=g^{\dagger}(g(a))$ or $g\big(g^{\dagger}(g(a))\big)-g(t)<\frac{\varepsilon}{2}$ and $s=g^{\dagger}(g(b))$ or $g\big(g^{\dagger}(g(b))\big)-g(s)<\frac{\varepsilon}{2}$, we conclude that $\|f\circ g-\chi_{[a,b)}\|_{L_g^1(\mathbb{R})}<\varepsilon$.

Since $\varepsilon$ was fixed arbitrarily, $\chi_{[a,b)}\in C$.

If $a\geqslant g^{\dagger}(g(b))$, either $g$ is constant on $[a,b]$, so it is enough to take $f=0$, or $a\in D_g$ and $g$ is constant on $(a,b]$ and it is enough to take $t$ as before and $f=P_A$ where
\[ A=\{(g(t),0),(g(a),1),(g(b),0)\},\]
\---see Figure~\ref{fig:proof:f:step1:ex3}.

\begin{figure}[h]
	\centering
	\resizebox{0.6\textwidth}{!}{%
		\begin{circuitikz}
			\tikzstyle{every node}=[font=\LARGE]
			\draw [ line width=2pt](-2.75,3) to[short] (15,3);
			\draw [line width=2pt, short] (-0.75,8.5) -- (-0.75,2.5);

			\draw [ color={rgb,255:red,17; green,113; blue,187}, line width=2pt, dashed] (4.4,3.2) -- (4.4,8.28);

			\node at (-0.75,3) [squarepole,line width=1.5pt] {};
			\node [font=\Large] at (-1,3.3) {$0$};
			\node at (-0.75,8.25) [squarepole,line width=1.5pt] {};
			\node [font=\Large] at (-1,8.2) {$1$};


			\node at (2,3) [squarepole,line width=1.5pt] {};
			\node [font=\LARGE] at (2.4,2.5) {$g(t)$};

			\node at (4.4,3) [squarepole,line width=1.5pt] {};
			\node at (4.4,8.25) [circ, color={rgb,255:red,17; green,113; blue,187},line width=0.001pt] {};
			\node [font=\LARGE] at (4.7,2.5) {$g(a)$};


			\node at (12.25,3) [squarepole,line width=1.5pt] {};
			\node at (12.25,3) [circ, color={rgb,255:red,17; green,113; blue,187},line width=0.001pt] {};
			\node [font=\LARGE] at (12.5,2.5) {$g(b)$};

			\draw [ color={rgb,255:red,17; green,113; blue,187}, line width=2.5pt] (-2.75,3) -- (2,3);
			\draw [ color={rgb,255:red,17; green,113; blue,187}, line width=2.5pt] (2,3) -- (4.4,8.25);
			\draw [ color={rgb,255:red,17; green,113; blue,187}, line width=2.5pt] (4.4,8.25) -- (12.25,3);
			\draw [ color={rgb,255:red,17; green,113; blue,187}, line width=2.5pt] (12.25,3) -- (15,3);


		\end{circuitikz}
	}%
	\caption{The continuous function $f$ in the case $a\geqslant g^{\dagger}(g(b))$ with $a\in D_g$, $t=g^\dagger(g(a))$, and $g\Big(g^{\dagger}(g(a))\Big)<g(a)$.}\label{fig:proof:f:step1:ex3}
\end{figure}
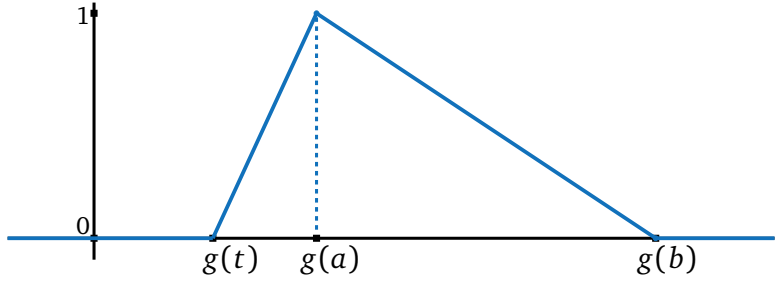

	 \emph{$\diamond$ Step 2: $\chi_{E}\in C$ with $E\in{\mathcal M}_g$ of finite $g$\--measure.} Fix $\varepsilon>0$. By Theorem~\ref{teoroutermeasureg}, and taking into account that $\mu_g(E)<\infty$, there exists a pairwise disjoint family of intervals $\{[a_n,b_n)\}_{n\in{\mathbb N}}$ such that
\[
E\subset\bigcup_{n\in{\mathbb N}}[a_n,b_n)\quad \text{and} \quad 0\leqslant\mu_g\big( \bigcup_{n\in{\mathbb N}}[a_n,b_n)\big) -\mu_g(E)<\frac{\varepsilon}{4}.
\]
Since $\mu_g\big( \bigcup_{n\in{\mathbb N}}[a_n,b_n)\big) <\infty$ and the family $\{[a_n,b_n)\}_{n\in{\mathbb N}}$ is pairwise disjoint, there exists $N\in{\mathbb N}$ such that $0\leqslant \mu_g\big( \bigcup_{n\in{\mathbb N}}[a_n,b_n)\big)-\mu_g\big( \bigcup_{n=1}^N[a_n,b_n)\big) <\frac{\varepsilon}{4}$, so
\[
 \big|\mu_g\big(\bigcup_{n=1}^N[a_n,b_n)\big) -\mu_g(E)\big| \leqslant \big|\mu_g\big( \bigcup_{n=1}^N[a_n,b_n)\big)- \mu_g\big(\bigcup_{n\in{\mathbb N}}[a_n,b_n)\big) \big| + \big| \mu_g\big( \bigcup_{n\in{\mathbb N}}[a_n,b_n)\big) -\mu_g(E) \big|<\frac{\varepsilon}{2}.
 \]
Now let us set $A:=\bigcup_{n=1}^{N}[a_n,b_n)$. Thus,
\[
E \setminus A \subset \bigcup_{n=N+1}^\infty [a_n,b_n), \text{ and }\mu_g(E \setminus A ) \leqslant \mu_g\big(\bigcup_{n=N+1}^\infty [a_n,b_n)\big)< \frac{\varepsilon}{4}.
\]
Also, as $A\setminus E \subset \big(\bigcup_{n\in{\mathbb N}}[a_n,b_n)\big)\setminus E$, and $E \subset \big(\bigcup_{n\in{\mathbb N}}[a_n,b_n)\big)$, it follows
\[
\mu_g(A\setminus E)\leqslant \mu_g\Big(\big(\bigcup_{n\in{\mathbb N}}[a_n,b_n)\big)\setminus E\Big)=\mu_g\big(\bigcup_{n\in{\mathbb N}}[a_n,b_n)\big)-\mu_g(E)<\frac{\varepsilon}{4}.
\]
Therefore,
\[
\mu_g(A\Delta E)=\mu_g(A\setminus E)+\mu_g(E\setminus A)<\frac{\varepsilon}{4}+\frac{\varepsilon}{4}=\frac{\varepsilon}{2}.
\]
Consequently,
\[
\|\chi_A-\chi_E \|_{L_g^1(\mathbb{R})}=\|\chi_{A\Delta E} \|_{L_g^1(\mathbb{R})}=\mu_g(A\Delta E)<\frac{\varepsilon}{2}.
\]
	 For every ${n\in{\mathbb N}}$, let $f_n$ be a continuous function such that $\|f_n\circ g-\chi_{[a_n,b_n)}\|_{L_g^1(\mathbb{R})}<\frac{\varepsilon}{2^{n+1}}$, as in Step~1. The function $f=\sum_{n=1}^N f_n$ is well\--defined, continuous and bounded. Furthermore,
	 \begin{align*}\|f\circ g\|_{L_g^1(\mathbb{R})}\leqslant \sum_{n=1}^N\|f_n\circ g\|_{L_g^1(\mathbb{R})}
\leqslant &\sum_{n=1}^N\|\chi_{[a_n,b_n)}\|_{L_g^1(\mathbb{R})}+\sum_{n=1}^N\|f_n\circ g-\chi_{[a_n,b_n)}\|_{L_g^1(\mathbb{R})}\\
\leqslant &\sum_{n=1}^N\mu_g([a_n,b_n))+\sum_{n=1}^N\frac{\varepsilon}{2^{n+1}}\\
< & \mu_g(E)+\frac{\varepsilon}{4}+\frac{\varepsilon}{2}<\infty,\end{align*}
	 so $f\circ g\in L^1_g(I,{\mathbb R})$.

Since the intervals $\{[a_n,b_n)\}_{n\in \mathbb{N}}$ are pairwise disjoint, we have that $\chi_{A}=\sum_{n=1}^N \chi_{[a_n,b_n)}$. Hence,
\[
 \|f\circ g-\chi_A\|_{L_g^1(\mathbb{R})}= \big\|\sum_{n=1}^N f_n\circ g-\chi_{[a_n,b_n)}\big\|_{L_g^1(\mathbb{R})}\leqslant \sum_{n=1}^N\|f_n\circ g-\chi_{[a_n,b_n)}\|_{L_g^1(\mathbb{R})} \leqslant \sum_{n=1}^N \frac{\varepsilon}{2^{n+1}}<\frac{\varepsilon}{2}.
\]
 Finally,
	 \[  \|f\circ g-\chi_E\|_{L_g^1(\mathbb{R})}\leqslant \|f\circ g-\chi_A\|_{L_g^1(\mathbb{R})} + \|\chi_A-\chi_E \|_{L_g^1(\mathbb{R})} <\frac{\varepsilon}{2}+\frac{\varepsilon}{2}=\varepsilon.\]

	 Since $\varepsilon$ was fixed arbitrarily, $\chi_{E}\in C$.

	\emph{$\diamond$ Step 3: $C=L^1_g(I,{\mathbb R})$.} Every $g$\--integrable simple function is a linear combination of the functions considered in Step 2. Therefore, given that $C$ is a vector space, they belong to $C$. Now, every function in $L^1_g(I,{\mathbb R})$ is the limit of simple functions. Since $C$ is closed, $C=L^1_g(I,{\mathbb R})$.

2. We prove the case where $I={\mathbb R}$ (the rest would be analogous). Let $f:I\to[c,d]$ be a $g$\--integrable function and fixed $\varepsilon>0$. Then, we can choose $\widetilde h\in{\mathcal B}{\mathcal C}_g({\mathbb R},{\mathbb R})$ such that $\|\widetilde h-f\|_{L_g^1(\mathbb{R})}< \varepsilon$. Defining $h=\min\{d,\max\{\widetilde h,c\}\}$ we get that $h\in{\mathcal B}{\mathcal C}_g({I},[c,d])$ and
	\[  \| h-f\|_{L_g^1(\mathbb{R})}\leqslant\|\widetilde h-f\|_{L_g^1(\mathbb{R})}< \varepsilon.\]

3. Let us fix $\varepsilon>0$, and define $p:=g^{\dagger}\circ g$. We will start studying the case where $g$ has a finite number of jumps $j\geqslant 0$ in $[a,b)$. The general case where the derivator~$g$ has infinitely many jumps will be studied via approximation by a derivator with finitely many jumps in $[a,b)$.

We set
\begin{equation}\label{eq:ell}
\ell:=\inf\{p^{(n)}(b): p^{(n)}(b)\geqslant a^*,\ n\in\mathbb{Z},\ n\geqslant 0\},
\end{equation}
where $p^{(n)}(b)=\overbrace{p \circ p \circ \dots \circ p}^{n\ \text{times}}(b)$, with $p^{(0)}(b)=b$. As $a^* \notin D_g$, and $g$ has a finite number of jumps, then $\ell \in (a^*,b]$, and $\ell$ is the first left accumulation point of~$g$ from the left of~$b$ in $(a^*,b]$ at which~$g$ is not constant on any left\--neighborhood of~$\ell$. Moreover, since $a^* \notin D_g$, then $a \notin D_g$ and $a^*$ is a right accumulation point of~$g$ at which~$g$ is not constant on any right\--neighborhood of~$a^*$. Thus, we can choose $s\in(a^*,\ell)$ and $r\in(a^*,s)$ such that
\[
0<g(r)-g(a)<\frac{\varepsilon}{3(d-c)},
\]
and
\[
0<g(\ell)-g(s)<\frac{\varepsilon}{3(d-c)}.
\]
Let us take $\widetilde h\in{\mathcal B}{\mathcal C}_g([r,s],[c,d])$ such that $\int_{{[r,s)}}|f|_{{[r,s)}}-\widetilde h|\operatorname{d} \mu_g<\frac{\varepsilon}{3}$, and define the sets:
\[
A_1:=\left\{\big(g(a),\alpha\big),\big(g(r),\widetilde h(r)\big)\right\},
\]
and
\[
A_2:= \begin{dcases}
 \left\{\big(g(s),\widetilde h(s)\big),\big(g(b),\beta\big)\right\}, & \mbox{if } [\ell,b)\cap D_g=\emptyset, \\
 \left\{\big(g(s),\widetilde h(s)\big)\right\}\cup\left\{\big(g(r_n),f(r_n)\big)\right\}_{n=1}^{m}\cup\left\{\big(g(b),\beta\big)\right\}, & \mbox{if } \{r_n\}_{n=1}^{m}:=[\ell,b)\cap D_g\ne\emptyset.
 \end{dcases}
\]
In the case where $\{r_n\}_{n=1}^{m}:=[\ell,b)\cap D_g\ne\emptyset$, the sequence $\{r_n\}_{n=1}^{m}$ is well\--ordered, increasing and satisfying $r_m=p(b)$, $r_{n}=p(r_{n+1})$ for $n=1,\dots,m-1$. Now, let us define the function $h:[a,b]\to\mathbb{R}$ by
	\[  h(t):=\begin{dcases}
(P_{A_1}\circ g)(t), & \mbox{{if} } t\in [a,r],\\
\widetilde h(t), &\mbox{{if} } t\in [r,s], \\
(P_{A_2}\circ g)(t), &\mbox{{if} } t\in [s,b].\end{dcases}
	\]
The function $h$ is, by construction, $g$\--continuous and $h([a,b])\subset[c,d]$. Furthermore,
\begin{align*}\int_{[a,b)}|f-h|\operatorname{d}\mu_g= & \int_{[a,r)}| P_{A_1}\circ g-f|_{[a,r)}|\operatorname{d}\mu_g+\int_{[r,s)}|\widetilde h-f|_{[r,s)}|\operatorname{d}\mu_g\\
 & \;\;+\int_{[s,\ell)}|P_{A_2}\circ g-f|_{[a,r)}|\operatorname{d}\mu_g+\int_{[\ell,b)}|P_{A_2}\circ g-f|_{[\ell,b)}|\operatorname{d}\mu_g\\
\leqslant & (d-c)\mu_g([a,r))+\frac{\varepsilon}{3}+(d-c)\mu_g([s,\ell))\\
= & (d-c)\big(g(r)-g(a)\big)+\frac{\varepsilon}{3}+ (d-c)\big(g(\ell))-g(s)\big)<\varepsilon.
\end{align*}

We now study the general case where $g$ has infinitely many jumps in $[a,b)$. First we claim that for a fixed $\eta>0$, there exists a left\--continuous nondecreasing~$G:{[a,b]}\to\mathbb{R}$ that has only finitely many discontinuities on $[a,b)$, satisfying
\[
\|\mu_g -\mu_G\|_{TV}={(\mu_g -\mu_G)([a,b))}<\eta,
\]
where $\mu_g$, $\mu_G$ are the Lebesgue\--Stieltjes measures generated by~$g$ and~$G$ respectively, and $\|\cdot\|_{TV}$ is the total variation norm on $\mathcal{M}([a,b],\mathcal{B}([a,b]))$ the Banach space of all signed measures of bounded variation defined on $\mathcal{B}([a,b])$ the $\sigma$\--algebra associated to the usual topology on $[a,b]$.

Let $\{r_n\}_{n=1}^\infty=D_g\subset [a, b)$ denote the countable set of discontinuities of $g$. At each discontinuity point $r_n \in D_g$, the jump of $g$ is given by
\[
\Delta^+ g(r_n) := g(r_n^+) - g(r_n)>0.
\]
Since $g$ is nondecreasing, the total variation of its jump discontinuities is finite:
\[
\sum_{n=1}^\infty \Delta^+ g(r_n) \leqslant g(b) - g(a) < \infty.
\]
Hence, for every $\eta >0$, there exists $N \in \mathbb{N}$ such that
\[
\sum_{n > N} \Delta^+ g(r_n) < \eta.
\]
Now, let us define the function $G: [a,b]\to \mathbb{R}$ by
\[
G(x) := g(x) - \sum_{\substack{n > N \\ r_n \in [a, x)}} \Delta^+ g(r_n), \quad \text{ for all } x \in [a,b] .
\]
Since $g$ is left\--continuous and nondecreasing and we are just subtracting some of the jumps when they happen, $G$~is also left\--continuous and nondecreasing\---see \cite{Fernandez2020}. Furthermore, $G$ has finitely many discontinuities, namely $\{r_n\}_{n=1}^N$.

 Let $\mu_G$ be the Lebesgue\--Stieltjes measure associated to~$G$. Then, the measure $\mu_g - \mu_G$ is given by
\[
\mu_g - \mu_G = \sum_{n > N} \Delta^+ g(r_n) \cdot \delta_{r_n},
\]
where $\delta_{r_n}$ denotes the Dirac mass at $r_n$. In particular, $\mu_g - \mu_G$ is a purely atomic measure supported on the set $\{r_n\}_{n>N} \subset [a,b)$, with total variation
\begin{equation}\label{eq:proof:density:<eta}
 \|\mu_g - \mu_G\|_{TV} = \sum_{\substack{n > N \\ r_n \in [a,b)}} \Delta^+ g(r_n) < \eta.
\end{equation}
Hence, the desired approximation holds.

In the sequel, we want to prove that there exists $h\in\mathcal{BC}_g([a,b],[c,d])$ such that $h(a)=\alpha$, $h(b)=\beta$ and $\|f-h\|_{L_g^1([a,b))}<\varepsilon$. First, let us fix a function~$G:[a,b]\to\mathbb{R}$ such that~\eqref{eq:proof:density:<eta} holds for
\[
\eta=\frac{\varepsilon}{2(d-c)}.
\]
We then consider the positive finite measure~$\nu:=\mu_g-\mu_G$, which is purely atomic and supported on the truncated set of jump points of~$g$. Observe that $f$ belongs to $L^{1}_{\nu}([a,b),\mathbb{R})$, the Banach space of integrable functions with respect to~$\nu$.

The function~$G$ is left\--continuous, nondecreasing, and $D_G\subset D_g$, thus $f$ is $G$\--integrable on $[a,b)$ as well. By the first part of the proof of point 2, there exists $h\in \mathcal{BC}_G([a,b],[c,d])$ such that $h(a)=\alpha$, $h(b)=\beta$ and
\[
\|f-h\|_{L_G^1([a,b))}<\frac{\varepsilon}{2}.
\]
As we have that $C_g\subset C_G$ and $D_G\subset D_g$, it follows by~\cite[Proposition~3.9]{MT} that~$G$ is $g$\--continuous. With this in mind, and since $h\in \mathcal{BC}_G([a,b],[c,d])$, we conclude that $h\in \mathcal{BC}_g([a,b],[c,d])$. Finally, since $\mu_g=\nu+\mu_G$, $|f-h|$ is bounded by $d-c$ and, taking~\eqref{eq:proof:density:<eta} into account for the $\eta$ defined above, we obtain
\begin{align*}
 \|f-h\|_{L_g^1([a,b))} = \int_{[a,b)} |f-h| \operatorname{d}\nu+ \int_{[a,b)} |f-h| \operatorname{d}\mu_G \leqslant & (d-c)\nu([a,b))+ \int_{[a,b)} |f-h| \operatorname{d}\mu_G\\
 <& (d-c) \|\mu_g-\mu_G\|_{TV}+\frac{\varepsilon}{2}<\varepsilon.
\end{align*}
Hence, we obtain $\|f-h\|_{L_g^1([a,b))}< \varepsilon$.

4. We proceed as in the previous point. Let $f:[a,b]\to[c,d]$ be a $g$\--integrable function and fix $\varepsilon>0$. Let us consider $\ell$ as in~\eqref{eq:ell}. Without loss of generality, let us assume that~$g$ has a finite number of jumps ${j}\geqslant 0$ in $[a,b)$. As $a^*\in D_g$, then $\ell\in[a^*,b]$. Thus, we distinguish two cases:

\textbf{Case 1: }$\ell=a^*$. Then, $g$ is a simple function with finitely many jumps $\{r_n\}_{n=1}^{j}\subset [a,b)$, ${j}\geqslant 1$ such that $r_k<r_{k+1}$ for $k=1,\dots, j-1$. Define the set
\[
A:=\left\{\big(g(a),f(a^*)\big)\right\} \cup \left\{\big(g(r_n),f(r_n)\big)\right\}_{n=1}^{j} \cup \left\{\big(g(b),\beta\big)\right\}.
\]

Now, let us define the function $h:[a,b]\to\mathbb{R}$ by $h:=P_A\circ g$. The function $h$ is, by construction, $g$\--continuous and $h([a,b])\subset[c,d]$. Furthermore,
\[
\int_{[a,b)}|f-h|\operatorname{d}\mu_g= \int_{[a^*,b)}| P_A\circ g-f|_{[a^*,b)}|\operatorname{d}\mu_g= \sum_{n=1}^{j}\int_{\{r_n\}}| P_A\circ g-f|_{[a^*,b)}|\operatorname{d}\mu_g = 0\leqslant \varepsilon.
\]

\textbf{Case 2:} $a^*<\ell$. In this case, we can take $s\in(a^*,\ell)$ such that
\[
0<g(\ell)-g(s)<\frac{\varepsilon}{2(d-c)}.
\]
Moreover, let us take $\widetilde h\in{\mathcal B}{\mathcal C}_g([a^*,s],[c,d])$ such that $\int_{{[a^*,s)}}|f|_{{[a^*,s)}}-\widetilde h|\operatorname{d} \mu_g<\frac{\varepsilon}{2}$, and define the set:
\[
A:= \begin{dcases}
 \left\{\big(g(s),\widetilde h(s)\big),\big(g(b),\beta\big)\right\}, & \mbox{if } [\ell,b)\cap D_g=\emptyset, \\
 \left\{\big(g(s),\widetilde h(s)\big)\right\} \cup\left\{\big(g(r_n),f(r_n)\big)\right\}_{n=1}^{m}\cup\left\{\big(g(b),\beta\big)\right\}, & \mbox{if } \{r_n\}_{n=1}^{m}:=[\ell,b)\cap D_g\ne\emptyset.
 \end{dcases}
\]
In the case where $\{r_n\}_{n=1}^{m}:=[\ell,b)\cap D_g\ne\emptyset$, the sequence $\{r_n\}_{n=1}^{m}$ is well\--ordered.
Now, let us define the function $h:[a,b]\to\mathbb{R}$ by
	\[  h(t):=\begin{dcases}
f(a^*), & \mbox{{if} } t\in[a,a^*],\\
\widetilde h(t), &\mbox{{if} } t\in (a^*,s],\\
P_A\circ g (t), &\mbox{{if} } t\in (s,b].\end{dcases}
	\]
The function $h$ is, by construction, $g$\--continuous (observe that $a^*\in D_g$) and $h([a,b])\subset[c,d]$. Furthermore,
\begin{align*}
\int_{[a,b)}|h-f|\operatorname{d}\mu_g= & \int_{[a,a^*)}\big|h-f|_{[a,a^*)}\big|\operatorname{d}\mu_g+ \int_{[a^*,s)}\big|\widetilde h-f|_{[a^*,s)}\big|\operatorname{d}\mu_g\\
&\quad+ \int_{[s,\ell)} \big|f|_{[s,\ell)}-P_A\circ g\big|\operatorname{d}\mu_g +\int_{[\ell,b)}\big|P_A\circ g-f|_{[\ell,b)}\big|\operatorname{d}\mu_g\\
< & \frac{\varepsilon}{2}+(d-c)\mu_g([s,\ell))=\frac{\varepsilon}{2}+(d-c)\big(g(\ell)-g(s)\big)<\varepsilon.
\end{align*}
To conclude the general case where $g$ has infinitely many jumps in $[a,b)$, one can apply the same approximation as in the proof of the previous point.
\end{proof}
\begin{rem}
	Although the approximation Theorem~\ref{density} was stated in terms of $g$\--continuous functions, note that the functions actually used are $g$\--uniformly continuous.
\end{rem}

In order to prove that condition~\eqref{eq:thm:varphi >0} in Theorem~\ref{thmfcf} is optimal, we need to prove the following lemma.
\begin{lem}\label{lem:lim f(x_n)-f(x_n+1)/h(x_n)-h(x_n+1)=0}
Let $t\in{\mathbb R}$, $A\subset(t,\infty)$ such that $t\in A'$, $B$ dense in $A$ and $f,h :A\to\mathbb{R}$ such that $h(s)\ne 0$ for every $s\in A$ and $f/h$ is left\--continuous on $A$. Assume that
 \[
 \lim_{s\to t^+}f(s)= \lim_{s\to t^+}h(s)=\liminf_{s\to t^+}\frac{f(s)}{h(s)}=0.
 \]
 Then, there exists a decreasing sequence $(x_n)_{n\in\mathbb{N}}\subset {B}$ such that $x_n\xrightarrow{n\to\infty} t^+$, and
\[
\lim_{n\to+\infty} \frac{f(x_n)-f(x_{n+1})}{h(x_n)-h(x_{n+1})}=0.
\]
Analogously, if $C\subset(-\infty,t)$ is such that $t\in C'$, $D$ is dense in $C$, $f,h :C\to\mathbb{R}$ are such that $h(s)\ne 0$ for every $s\in B$, $f/h$ is left continuous and
\[
\lim_{s\to t^-}f(s)= \lim_{s\to t^-}h(s)=\liminf_{s\to t^-}\frac{f(s)}{h(s)}=0,
 \]
 then there exists an increasing sequence $(x_n)_{n\in\mathbb{N}}\subset D$ such that $x_n\xrightarrow{n\to\infty} t^-$, and
\[
\lim_{n\to+\infty} \frac{f(x_n)-f(x_{n+1})}{h(x_n)-h(x_{n+1})}=0.
\]
\end{lem}
\begin{proof}
We prove the first part; the second follows analogously. Since $\liminf\limits_{s\to t^+}\frac{f(s)}{h(s)}=0$, there exist a decreasing sequence $(\widetilde y_n)_{n\in\mathbb{N}}$ in $A$ such that
 $\widetilde y_n\xrightarrow{n\to\infty} t^+$, and
\[
\lim_{n\to+\infty}\frac{f(\widetilde y_n)}{h(\widetilde y_n)}=\liminf_{s\to t^+}\frac{f(s)}{h(s)}=0.
\]
Since $f/h$ is left\--continuous, for every $n\in{\mathbb N}$ there exists $z_n\in(t,\widetilde y_n)$ such that $|f(t_n)/h(t_n)-f(\widetilde y_n)/h(\widetilde y_n)|<2^{-n}$ for every $t_n\in(z_n,\widetilde y_n)$. Since $B$ is dense in $A$, we can take $ y_n\in (z_n,\widetilde y_n)\cap B\ne\emptyset$ for every ${n\in{\mathbb N}}$. Given that $f(\widetilde y_n)/h(\widetilde y_n)\to 0$, it is clear that $f(y_n)/h(y_n)\to 0$.

Given that $h(s)\neq 0$ for all $s\in A$, observe that the function

\begin{center}
	\begin{tikzcd}[
		,row sep = 0ex
		,/tikz/column 1/.append style={anchor=base east}
		,/tikz/column 2/.append style={anchor=base west}
		]
	\mathbb{R} \times \left( \mathbb{R} \setminus \{h(s)\}\right) \arrow[r, "H_{s,c}"] & \mathbb{R} \\
	(u,v) \arrow[mapsto, r] & \dfrac{c-u}{h(s)-v}
		\end{tikzcd}
\end{center}
 is continuous at the point $(0,0)$ for any $(s,c)\in A\times\mathbb{R}$ fixed.

Let $x_1\in B$ be fixed. As $\lim_{s\to t^+}f(s)=\lim_{s\to t^+}h(s)=0$, and $H_{x_1,f(x_1)}$ is continuous at $(0,0)$, we can choose $x_2\in \{y_n\}_{n\in\mathbb{N}}$ satisfying $x_2<\min\{x_1,t+1\}$ such that $(f(x_2),h(x_2))$ are sufficiently close to $(0,0)$, to guarantee that $|H_{x_1,f(x_1)}(f(x_2),h(x_2))-H_{x_1,f(x_1)}(0,0)|<1$, that is,
\[
\left|\frac{f(x_1)-f(x_2)}{h(x_1)-h(x_2)} - \frac{f(x_1)}{h(x_1)}\right|<1.
\]
Next, we choose similarly $x_3\in \{y_n\}_{n\in\mathbb{N}}$ such that $x_3<\min\left\{x_2,t+\frac{1}{2}\right\}$, and
\[
\left|\frac{f(x_2)-f(x_3)}{h(x_2)-h(x_3)} - \frac{f(x_2)}{h(x_2)}\right|<\frac{1}{2}.
\]
We repeat this process inductively: for each $n\in\mathbb{N}$, given the continuity of $H_{x_n,f(x_n)}$ at $(0,0)$, choose $x_{n+1} \in \{y_n\}_{n\in\mathbb{N}}$ such that $x_{n+1}<\min\{x_n,t+\frac{1}{n}\}$, and
\[
\left|\frac{f(x_n)-f(x_{n+1})}{h(x_n)-h(x_{n+1})} - \frac{f(x_n)}{h(x_n)}\right|<\frac{1}{n}.
\]
By construction, $(x_n)_{n\in\mathbb{N}}$ is a decreasing subsequence of $(y_n)_{n\in\mathbb{N}}$. We have $x_n\xrightarrow{n\to\infty} t^+$, and hence we obtain
\[
\lim_{n\to+\infty}\frac{f(x_n)}{h(x_n)}=0.
\]
Consequently, we conclude that
\[
\lim_{n\to+\infty}\frac{f(x_n)-f(x_{n+1})}{h(x_n)-h(x_{n+1})}=0.\qedhere
\]
\end{proof}

Now, in the next theorem, we prove that Condition~\eqref{eq:thm:varphi >0} in Theorem~\ref{thmfcf} is optimal.
\begin{thm}\label{thmnecessarycondition} Let $g\in \operatorname{BV}^{\operatorname{lc}}([a,b],{\mathbb R})$ be such that $b\notin N_g^+$.
Consider $\varphi$ defined as in~\eqref{eq:thm:varphi function}. If $\varphi(t)=0$ for some $t\in[a,b]$, then there exists $f \in \mathcal{BC}_g([a,b],\mathbb{R})$ such that the function $F(s):=\int_{[a,s)} f\, \operatorname{d} \mu_g \in \mathbb{R}$, $s\in[a,b]$, is not $g$\--differentiable at~$t$.
\end{thm}
\begin{proof}Let $t\in[a,b]$. By Remark~\ref{rem:phi(t)>0 is required on continuity points}, we have ${t^*}\notin D_g$, thus, $\lim\limits_{s\to {t^*}^\pm}|g(s)-g({t^*})|=\lim\limits_{s\to {t^*}^\pm}\widetilde{g}(s)-\widetilde{g}({t^*})=0$. Moreover, since $\widetilde{g}$ is nondecreasing, it follows that $\widetilde{g}(s)-\widetilde{g}({t^*})\ne 0$ for all $s>{t^*}$ with ${t^*\in[a,b]\setminus (D_g\cup N_g^-)}$, and all $s<{t^*}$ with ${t^*\in[a,b]\setminus (D_g\cup N_g^+)}$. By the definition of $\varphi$ in~\eqref{eq:thm:varphi function} and since $\varphi(t)=0$, {we have that $\varphi(t^*)=0$}, and, given that $[a,b]\setminus D_g$ is dense in $[a,b]$, it follows from Lemma~\ref{lem:lim f(x_n)-f(x_n+1)/h(x_n)-h(x_n+1)=0} that there exists a strictly decreasing or increasing sequence there exists a decreasing sequence $(x_n)_{n\in\mathbb{N}}\subset {(t^*,b]\backslash D_g}$ such that $x_n\xrightarrow{n\to\infty} (t^*)^+$ and
\begin{equation}\label{eq:proof:direct appli. of lim f(x_n)-f(x_n+1)/h(x_n)-h(x_n+1)=0}
 \lim_{n\to+\infty} \frac{| g(x_n)-g({t^*})|-|g(x_{n+1})-g({t^*})|}{\widetilde g(x_n)- \widetilde g(x_{n+1})}=0,
\end{equation}
Observe that this implies that $|g(x_n)-g({t^*})|\to 0$. Without loss of generality, we will assume that $(g(x_n))_{n\in{\mathbb N}}$ is strictly monotone, $(x_n)_{n\in{\mathbb N}}$ is decreasing (we can guarantee this by passing to a subsequence if necessary), and that $x_1=b$. In particular, since $(g(x_n))_{n\in{\mathbb N}}$ is strictly monotone and converges to $g({t^*})$, we obtain that $(|g(x_n)-g({t^*})|)_{n\in \mathbb{N}}$ tends to~$0$ from the right, and also
\begin{align*}
 |g(x_n)-g({t^*})|-|g(x_{n+1})-g({t^*})|= &\begin{dcases}
 g(x_n)-g(x_{n+1}), & \mbox{if $(g(x_n))_{n\in{\mathbb N}}$ is decreasing,} \\
 g(x_{n+1})- g(x_n), & \mbox{if $(g(x_n))_{n\in{\mathbb N}}$ is increasing,}
 \end{dcases}\\
 = & |g(x_n)-g(x_{n+1})|.
\end{align*}
Thus,~\eqref{eq:proof:direct appli. of lim f(x_n)-f(x_n+1)/h(x_n)-h(x_n+1)=0} becomes
\begin{equation*}
 \lim_{n\to+\infty} \frac{|g(x_n) - g(x_{n+1})|}{\widetilde g(x_n) - \widetilde g(x_{n+1})}=0.
\end{equation*}

Let us set, for ${n\in{\mathbb N}}$,
 \[  M_n:=\sqrt{\frac{|g(x_{n+1})-g(x_n)|}{\widetilde g(x_n)-\widetilde g(x_{n+1})}}>0.\]
Observe that $M_n\to 0$ and define
\[
\varepsilon_n:=|g(x_{n+1})-g(x_n)|>0,
\]
for every $n\in{\mathbb N}$. Observe that $\varepsilon_n\to 0$ and that
 \[  \frac{\sum_{k=n}^\infty\varepsilon_k}{|g(x_n)-g({t^*})|}=\frac{|g(x_n)-g(t^*)|}{|g(x_n)-g({t^*})|}=1\quad \text{for $n\in{\mathbb N}$.}\]

For each $n\in{\mathbb N}$, let us denote $g_n:=g|_{[x_{n+1},x_n]}$ and $\widetilde{g_n}:=\widetilde{g}|_{[x_{n+1},x_n]}$. Observe that, since $x_k\not\in D_g$ for any $k\in{\mathbb N}$, $\mu_{\widetilde{g_n}}=|\mu_{g_n}|$, and $\mu_{g_n}=\mu_g|_{[x_{n+1},x_n)}$. Now, by applying Theorem~\ref{density}, points 3 and 4, on $\widetilde{g_n}$, we can construct a bounded $g_n$\--continuous function $u_n:[x_{n+1},x_n]\to[0,1]$ such that
\[
 u_{n}(x_{n+1})= \begin{dcases}
 0, &\mbox{if } x_{n+1}^*\notin D_g,\\
\chi_{A_g^+}(x_{n+1}^*), &\mbox{if } x_{n+1}^*\in D_g,
 \end{dcases}
 \quad
 u_n(x_n)=0,
 \quad\text{and} \quad
 M_n\int_{[x_{n+1},x_n)}|u_n-\chi_{A_g^+}|\operatorname{d} |\mu_{g_n}|<\frac{\varepsilon_n}{2}.
\]
Let
\[
f_1(x):=\begin{dcases}
0, &\mbox{if } x\in[a,t^*],\\
M_nu_n(x), &\mbox{if } x\in{(x_{n+1},x_n]},\ {n\in{\mathbb N}}.
\end{dcases}
\]
Observe that $f_1$ is a well\--defined $g$\--continuous function on $[a,b]$. Indeed, if $t^*\neq a$, $f_1\equiv 0$ on $[a,t^*]$. Thus, for each $x\in [a,t^*)$, and every $\epsilon>0$ fixed, there exists $\delta=\widetilde g(t^*)-\widetilde g(x)>0$ such that for all $s\in [a,b]$
\[
|\widetilde g(s)-\widetilde g(x)|<\delta \Rightarrow |f_1(s) -f_1(x)|=0< \epsilon,
\]
and therefore, $f_1$ is $g$\--continuous at each $x\in[a,t^*)$. Moreover, for each $n\in\mathbb{N}$, $u_n$ is $g_n$\--continuous and, $u_n(x_{n+1})=u_{n+1}(x_{n+1})$ whenever $x_{n+1}^*\notin D_g$. This implies that $f_1$ is $g$\--continuous from the left and from the right at $x_{n+1}$ for each $n$, which guarantees $g$\--continuity at the gluing point $x_{n+1}$, and we conclude that $f_1$ is continuous on $(t^*,b]$.

Finally, let us check that $f_1$ is $g$\--continuous at $t^*$. Let $\varepsilon>0$ be fixed. Since $M_n\to 0$, we choose $N\in\mathbb N$ such that $M_n<\varepsilon$ for all $n\geqslant N$. Since $x_N>t^*$, define $\delta:=\widetilde g(x_N)-\widetilde g(t^*)>0$, and let $s\in[a,b]$ satisfy $|\widetilde g(s)-\widetilde g(t^*)|<\delta$. If $s\leqslant t^*$, then $f(s)=f(t^*)=0$. Assume now that $s>t^*$. Thus, there exists $n\in\mathbb N$ such that $s\in(x_{n+1},x_n]$. We claim that $n\geqslant N$. Indeed, if $n<N$, then $x_{n+1}\geqslant x_N$ because $(x_n)_{n\in\mathbb N}$ is decreasing. Since $s>x_{n+1}$, and $\widetilde g$ is nondecreasing,
\[
\widetilde g(s)\geqslant \widetilde g(x_{n+1})\geqslant \widetilde g(x_N).
\]
Thus, $\widetilde g(s)-\widetilde g(t^*)\geqslant \widetilde g(x_N)-\widetilde g(t^*)=\delta$ which contradicts the choice of~$s$. Hence $n\geqslant N$, and we obtain
\[
|f_1(s)-f_1(t^*)|
\leqslant M_n|u_n(s)|\leqslant M_n\leqslant M_N<\varepsilon.
\]

 Analogously, let us define a bounded $g_n$\--continuous function $v_n:[x_{n+1},x_n]\to[-1,0]$ such that
\[
 v_{n}(x_{n+1})= \begin{dcases}
 0, &\mbox{if } x_{n+1}^*\notin D_g,\\
-\chi_{A_g^-}(x_{n+1}^*), &\mbox{if } x_{n+1}^*\in D_g,
 \end{dcases}
 \quad
 v_n(x_n)=0,
 \quad\text{and} \quad
 M_n\int_{[x_{n+1},x_n)}|v_n+\chi_{A_g^-}|\operatorname{d} |\mu_{g_n}|<\frac{\varepsilon_n}{2}.
\]
Let
\[
f_2(x):=\begin{dcases}
0, &\mbox{if } x\in[a,t^*],\\
M_nv_n(x), &\mbox{if } x\in(x_{n+1},x_n],\ n\in{\mathbb N}.
\end{dcases}
\]

Let $f:=f_1+f_2$, $F(s):=\int_{[a,s)}f\operatorname{d}\mu_g$. Observe that $f_1$ is positive on~$A_g^+$, and $f_2$ is negative on~$A_g^-$, therefore, $F$ is nondecreasing. Moreover, $f$ is a bounded $g$\--continuous function and, for ${n\in{\mathbb N}}$, since $x_n\notin D_g$, then
\begin{align*}& \left|F(x_n)-F(x_{n+1})-M_n(\widetilde g(x_n)-\widetilde g(x_{n+1}))\right|
\\=& \left|\int_{[x_{n+1},x_n)}f \operatorname{d} \mu_g- M_n|\mu_{g_n}|([x_{n+1},x_n))\right|
\\= & \left|\int_{[x_{n+1},x_n)}f \operatorname{d} \mu_{g_n}- M_n\mu_{g_n}^+((x_{n+1},x_n))-M_n\mu_{g_n}^-((x_{n+1},x_n))\right|
\\= & \left|\int_{[x_{n+1},x_n)}M_n(u_n+v_n) \operatorname{d} \mu_{g_n}- M_n \int_{[x_{n+1},x_n)}\chi_{A_g^+}\operatorname{d} \mu_{g_n} + M_n \int_{[x_{n+1},x_n)}\chi_{A_g^-}\operatorname{d} \mu_{g_n}\right|
\\= & M_n\left|\int_{[x_{n+1},x_n)}(u_n+v_n-\chi_{A_g^+}+\chi_{A_g^-})\operatorname{d} \mu_{g_n}\right|\\ \leqslant & M_n\int_{[x_{n+1},x_n)}|u_n-\chi_{A_g^+}|\operatorname{d} |\mu_{g_n}|+M_n\int_{[x_{n+1},x_n)}|v_n+\chi_{A_g^-}| \operatorname{d} |\mu_{g_n}|<\varepsilon_n.\end{align*}
Thus, using the triangle inequality and the fact that $F$ and $\widetilde g$ are nondecreasing,
\[  M_n[\widetilde g(x_n)-\widetilde g(x_{n+1})]-\varepsilon_n=M_n|\widetilde g(x_n)-\widetilde g(x_{n+1})|-\varepsilon_n\leqslant|F(x_n)-F(x_{n+1})|=F(x_n)-F(x_{n+1}).\]
Hence,
\begin{align*}\left|\frac{F(x_n)-F({t^*})}{g(x_n)-g({t^*})}\right|= & \frac{F(x_n)-F({t^*})}{|g(x_n)-g({t^*})|}=\frac{\sum_{k=n}^\infty( F(x_k)-F(x_{k+1}))}{|g(x_n)-g({t^*})|}\\\geqslant & \frac{\sum_{k=n}^\infty (M_k [\widetilde g(x_k)-\widetilde g(x_{k+1})]-\varepsilon_k)}{|g(x_n)-g({t^*})|}\\= & \frac{\sum_{k=n}^\infty M_k [\widetilde g(x_k)-\widetilde g(x_{k+1})]}{|g(x_n)-g({t^*})|}- \frac{\sum_{k=n}^\infty\varepsilon_k}{|g(x_n)-g({t^*})|}=\frac{\sum_{k=n}^\infty M_k [\widetilde g(x_k)-\widetilde g(x_{k+1})]}{|g(x_n)-g({t^*})|}-1.
\end{align*}
Now, since $(g(x_n))_{n\in\mathbb{N}}$ is strictly monotone and converges to $g(t)$, either $|g(x_n)-g({t^*})|=g(x_n)-g({t^*})$ for every ${n\in{\mathbb N}}$ or $|g(x_n)-g({t^*})|=-[g(x_n)-g({t^*})]$ for every ${n\in{\mathbb N}}$. We consider the first case, as the second is analogous. In this case, we have that $(g(x_n)-g({t^*}))_{n\in{\mathbb N}}$ is a strictly decreasing function converging to zero and $\left( \sum_{k=n}^\infty M_k [\widetilde g(x_k)-\widetilde g(x_{k+1})]\right) _{n\in{\mathbb N}}$ converges to zero as well, so we can use the Stolz\--Cesàro Theorem \cite[Theorem 2.7.1]{choudary2014real} to deduce that
\begin{align*}
	\lim_{n\to\infty}\frac{\sum_{k=n}^\infty M_k [\widetilde g(x_k)-\widetilde g(x_{k+1})]}{g(x_n)-g({t^*})}= & \lim_{n\to\infty}\frac{\sum_{k=n+1}^\infty M_k [\widetilde g(x_k)-\widetilde g(x_{k+1})]-\sum_{k=n}^\infty M_k [\widetilde g(x_k)-\widetilde g(x_{k+1})]}{[g(x_{n+1})-g({t^*})]-[g(x_n)-g({t^*})]}\\ = & -\lim_{n\to\infty}\frac{M_n[\widetilde g(x_n)-\widetilde g(x_{n+1})]}{g(x_{n+1})-g(x_n)}=\lim_{n\to\infty}\frac{1}{M_n}=\infty.
	\end{align*}
Thus, $F$ cannot be $g$\--differentiable at $t$.
\end{proof}
\section*{Acknowledgments}
Lamiae Maia would like to express her sincere gratitude to Professor Fernando Adrián Fernández Tojo and to the Departamento de Estatística, Análise Matemática e Optimización, Facultade de Matemáticas, Universidade de Santiago de Compostela, for their warm hospitality during her research stay. Lamiae Maia also acknowledges the GI\--1561 research group “Nonlinear Differential Equations” at the University of Santiago de Compostela for funding this stay, during which the present article was initiated.
\section*{Funding}
Lamiae Maia was partially supported by the “National Center for Scientific and Technical Research (CNRST)” under the Excellence Research Fellowship Program, Grant No. 60UM5R2021, Morocco.

F. Adrián F. Tojo was supported by grant PID2020\--113275GB\--I00 funded by Xunta de Galicia, Spain, project ED431C 2023/12; and by MCIN/AEI/10.13039/ 501100011033, Spain, and by “ERDF A way of making Europe” of the “European Union”.

\bibliography{FullBibGD}

\begin{thebibliography}{10}
\providecommand{\url}[1]{{#1}}
\providecommand{\urlprefix}{URL }
\expandafter\ifx\csname urlstyle\endcsname\relax
  \providecommand{\doi}[1]{DOI~\discretionary{}{}{}#1}\else
  \providecommand{\doi}{DOI~\discretionary{}{}{}\begingroup
  \urlstyle{rm}\Url}\fi

\bibitem{AFNT}
Area, I., Fern{\'a}ndez, F.J., Nieto, J.J., F.~Tojo, F.A.: \emph{Concept and
  solution of digital twin based on a Stieltjes differential equation}.
\newblock Math. Methods Appl. Sci. \textbf{45}, 7451 -- 7465 (2022)

\bibitem{AL}
Athreya, K.B., Lahiri, S.N.: Measure Theory and Probability Theory.
\newblock Springer Science \& Business Media, New York (2006)

\bibitem{AP}
Aversa, V., Preiss, D.: \emph{Lusin's theorem for derivatives with respect to a
  continuous function}.
\newblock Proc. Amer. Math. Soc. \textbf{127}(11), 3229--3235 (1999)

\bibitem{BenedettoCzaja2009integration}
Benedetto, J.J., Czaja, W.: Integration and modern analysis, vol.~15.
\newblock Springer (2009)

\bibitem{Cacc}
Caccioppoli, R.: \emph{Sul lemma fondamentale del calcolo integrale}.
\newblock Atti Mere. Accad. Sci. Padova \textbf{50}, 93--98 (1934)

\bibitem{choudary2014real}
Choudary, A.D.R., Niculescu, C.P.: Real analysis on intervals.
\newblock Springer (2014)

\bibitem{D1918}
Daniell, P.J.: \emph{Differentiation with respect to a function of limited
  variation}.
\newblock Trans. Amer. Math. Soc. \textbf{19}(4), 353--362 (1918)

\bibitem{D1920}
Daniell, P.J.: \emph{{Stieltjes' Derivatives}}.
\newblock Proc. London Math. Soc. (3) \textbf{s2-30}(1), 187--198 (1930)

\bibitem{Feller}
Feller, W.: \emph{On differential operators and boundary conditions}.
\newblock Commun. Pure Appl. Math. \textbf{8}(1), 203--216 (1955)

\bibitem{FMarTo-OnFirstandSec}
Fernández, F.J., Márquez~Albés, I., Tojo, F.A.F.: \emph{On first and second
  order linear Stieltjes differential equations}.
\newblock J. Math. Anal. Appl. \textbf{511}(1), 126,010 (2022)

\bibitem{FernándezAlbésTojo2024FirstNSec}
Fernández, F.J., Márquez~Albés, I., Tojo, F.A.F.: \emph{On second-order
  linear Stieltjes differential equations with non-constant coefficients}.
\newblock Open Math. \textbf{22}(1), 20240,018 (2024)

\bibitem{FMT2025kernel-Stieltjes-deriv-space}
Fernández, F.J., Márquez~Albés, I., Tojo, F.A.F., Villanueva~Mariz, C.:
  \emph{On the kernel of the Stieltjes derivative and the space of bounded
  Stieltjes-differentiable functions}.
\newblock Electron. J. Qual. Theory Differ. Equ. (36), 1--41 (2025)

\bibitem{Fernandez2020}
Fernández, F.J., Tojo, F.A.F.: \emph{Numerical Solution of {Stieltjes}
  Differential Equations}.
\newblock Mathematics \textbf{8}(9) (2020)

\bibitem{FernandezTojoStieltjesBochnerSpaces}
Fernández, F.J., Tojo, F.A.F.: \emph{Stieltjes Bochner spaces and applications
  to the study of parabolic equations}.
\newblock J. Math. Anal. Appl. \textbf{488}(2), 124,079 (2020)

\bibitem{FP}
Frigon, M., López~Pouso, R.: \emph{Theory and applications of first-order
  systems of Stieltjes differential equations}.
\newblock Adv. Nonlinear Anal. \textbf{6}(1), 13--36 (2017)

\bibitem{FrTo}
Frigon, M., Tojo, F.A.F.: \emph{Stieltjes differential systems with
  nonmonotonic derivators}.
\newblock Bound. Value Probl. \textbf{2020}, 1--24 (2020)

\bibitem{Garg1992relativization}
Garg, K.M.: Relativization of some aspects of the theory of functions of
  bounded variation.
\newblock Polska Akademia Nauk, Instytut Matematyczny (1992)

\bibitem{Gradin}
Gradinaru, M.: \emph{On the derivative with respect to a function with
  applications to Riemann-Stieltjes integral}.
\newblock In: Collection: Seminar on Mathematical Analysis (Cluj-Napoca,
  1989-1990), pp. 21--28 (1989)

\bibitem{ThesisFLariviere}
Larivi{\`e}re, F.: \emph{Sur les solutions d’équations différentielles de
  Stieltjes du premier et du deuxième ordre}.
\newblock Mémoire de maîtrise, Université de Montréal (2018)

\bibitem{Liberman}
Liberman, J.: \emph{Théorème de Denjoy sur la dérivée d’une fonction
  arbitraire par rapport à une fonction continue}.
\newblock Recueil Mathématique \textbf{9(51)}, 221--236 (1941)

\bibitem{PMM}
L{\'o}pez~Pouso, R., Márquez~Albés, I., Monteiro, G.A.: \emph{Extremal
  solutions of systems of measure differential equations and applications in
  the study of Stieltjes differential problems.}
\newblock Electron. J. Qual. Theo.  (2018)

\bibitem{PR}
L{\'o}pez~Pouso, R., Rodr{\'i}guez, A.: \emph{{A New Unification of Continuous,
  Discrete and Impulsive Calculus through Stieltjes Derivatives}}.
\newblock Real Anal. Exchange \textbf{40}(2), 319 -- 354 (2015)

\bibitem{PM}
{López Pouso}, R., {Márquez Albés}, I.: \emph{General existence principles
  for Stieltjes differential equations with applications to mathematical
  biology}.
\newblock J. Differ. Equations \textbf{264}(8), 5388--5407 (2018)

\bibitem{PM2}
{López Pouso}, R., Márquez~Albés, I.: \emph{Resolution methods for
  mathematical models based on differential equations with Stieltjes
  derivatives}.
\newblock Electron. J. Qual. Theo. \textbf{2019}(72), 1--15 (2019)

\bibitem{PM3}
{López Pouso}, R., {Márquez Albés}, I.: \emph{Systems of Stieltjes
  differential equations with several derivators}.
\newblock Mediterr. J. Math. \textbf{16}(2), 51 (2019)

\bibitem{MEF1}
Maia, L., El~Khattabi, N., Frigon, M.: \emph{Existence and multiplicity results
  for first-order Stieltjes differential equations}.
\newblock Adv. Nonlinear Stud. \textbf{22}(1), 684--710 (2022)

\bibitem{MEF2}
Maia, L., El~Khattabi, N., Frigon, M.: \emph{Systems of Stieltjes differential
  equations and application to a predator-prey model of an exploited fishery}.
\newblock Discrete Cont. Dyn-a. \textbf{43}(12), 4244--4271 (2023)

\bibitem{MEF3}
Maia, L., El~Khattabi, N., Frigon, M.: \emph{Prolongation of solutions and
  Lyapunov stability for Stieltjes dynamical systems}.
\newblock Electron. J. Qual. Theory Differ. Equ. (19), 1--37 (2025)

\bibitem{ThesisMarquezAlbes}
M\'arquez~Albés, I.: \emph{Differential problems with Stieltjes derivatives
  and applications}.
\newblock Ph.D. thesis, Universidade de Santiago de Compostela (2021)

\bibitem{AlbesSlavikLogisticEq}
M\'arquez~Albés, I., Slavík, A.: \emph{The logistic equation in the context
  of Stieltjes differential and integral equations}.
\newblock Electron. J. Qual. Theory Differ. Equ \textbf{2023}, 1--35 (2023)

\bibitem{MT}
M\'arquez~Albés, I., Tojo, F.A.F.: \emph{Existence and uniqueness of solution
  for Stieltjes differential equations with several derivators}.
\newblock Mediterr. J. Math. \textbf{18}, 1--31 (2021)

\bibitem{MonteiroBianca-GenerDeriv2017}
Monteiro, G.A., Satco, B.: \emph{Distributional, differential and integral
  problems: equivalence and existence results}.
\newblock Electron. J. Qual. Theory Differ. Equ pp. 1--26 (2017)

\bibitem{Petrov}
Petrovsky, J.: \emph{Sur l'unicité de la fonction primitive par rapport à une
  fonction continue arbitraire}.
\newblock Rec. Math. Soc. Math. Moscou \textbf{41}(1), 48--59 (1934)

\bibitem{Tojo2025connectingSDEnODE}
Tojo, F.A.F.: \emph{On the connection between Stieltjes differential equations
  and ordinary differential equations}.
\newblock J. Math. Anal. Appl. \textbf{546}(1), 129,248 (2025)

\bibitem{Young}
Young, W.: \emph{On integrals and derivates with respect to a function}.
\newblock Proc. London Math. Soc. (3) \textbf{2}(1), 35--63 (1917)

\end{thebibliography}
\bibliographystyle{spmpsciper}

\end{document}